\documentclass[letterpaper,12pt]{article}

\usepackage{fancyhdr}
\usepackage{amsmath}
\usepackage{amsbsy}
\usepackage{amssymb}
\usepackage{amsthm}
\usepackage{amscd}
\usepackage{amsfonts}
\usepackage{supertabular}
\usepackage{graphics}
\usepackage{graphicx}
\usepackage{verbatim}
\usepackage{subfigure}
\usepackage{epsfig}
\usepackage{xspace}
\usepackage{euscript}
\usepackage{alltt}
\usepackage{boxedminipage}
\usepackage{float}
\usepackage{times}
\usepackage[colorlinks]{hyperref}
\usepackage[square,numbers]{natbib}
\usepackage{setspace}
\usepackage{bm}
\usepackage{bbm}
\usepackage[numbered,framed]{matlab-prettifier}
\usepackage[final]{pdfpages}
\usepackage{algorithm}
\usepackage{algorithmic}
\usepackage{authblk}
\usepackage{wasysym}

\begin{document}

\title{Explicit variational time integrator based on \\ rescaled Rodrigues parameters}
\author[1]{Caroline Baker}
\author[1,2]{Marcial Gonzalez  \thanks{Corresponding author: Marcial Gonzalez. E-mail: marcial-gonzalez@purdue.edu}}
\affil[1]{School of Mechanical Engineering, Purdue University, \break West Lafayette, IN 47907, USA}
\affil[2]{Ray W. Herrick Laboratories, Purdue University, \break West Lafayette, IN 47907, USA}

\maketitle
\begin{abstract}
	
We develop an explicit, second-order, variational time integrator for full body dynamics that preserves the momenta of the continuous dynamics, such as linear and angular momenta, and exhibits near-conservation of total energy over exponentially long times. In order to achieve these properties, we parametrize the space of rotations using exponential local coordinates represented by a rescaled form of the Rodrigues rotation vector and we systematically derive the time integrator from a  discrete Lagrangian function that yields discrete Euler-Lagrange equations amenable to explicit, closed-form solutions. By restricting attention to spherical bodies and Lagrangian functions with a quadratic kinetic energy and potential energies that solely depend on positions and attitudes, we show that the discrete Lagrangian map exhibits the same mathematical structure, up to terms of second order, of explicit Newmark or velocity Verlet algorithms, both known to be variational time integrators. These preserving properties, together with linear convergence of trajectories and quadratic convergence of total energy, are born out by two examples, namely the dynamics on $\mathrm{SO}(3)$ of a three-dimensional pendulum, and the nonlinear dynamics on $\mathrm{SE}(3)^n$ that results from the impact of a particle-binder torus against a rigid wall.	
	
\end{abstract}

\section{Introduction}

The formulation of explicit time-integration schemes for full body dynamics with exact conservation properties has been a longstanding goal of numerical analysis and computational mechanics. The reasons are two-fold. Firstly, integrators with a fixed time step cannot simultaneously preserve energy, the symplectic structure and other conserved quantities, such as linear and angular momenta \cite{GeMarsden1988LiePoissonIntegrators}. This restriction was overcome for dynamics on $\mathrm{E}(3)^n$ by Gonzalez and Ortiz \cite{EnergySteppingIntegrators} who proposed a method of approximation for Lagrangian mechanics consisting of replacing the Lagrangian function of the system by a sequence of approximate Lagrangians that can be solved exactly. Since these approximate solutions are the exact trajectories of the approximate Lagrangian, they are symplectic and preserve total energy and all momentum maps whose associated symmetries are preserved by the approximate Lagrangians. Secondly, the configuration manifold $\mathrm{SE}(3)$ is a Lie group and its representation is challenging since the space of finite three-dimensional rotations is not linear. A common approach for this representation is to introduce a linear approximation, at least locally, and use local coordinates in the form of, among others: (i) three angle coordinates, such as Euler angles, or (ii) exponential coordinates given by a rotation vector $\theta \in \mathbb{R}^3$ such that $R=\mathrm{e}^{S(\theta)}$. These descriptions often involve complicated trigonometric or transcendental expressions and introduce complexity in the analysis and computations.

The design of symplectic-momentum methods can be accomplished in a systematic and natural manner by recourse to a discrete version of Hamilton’s variational principle  based on an approximating discrete Lagrangian (cf., for example \cite{Marsden2001,Lew2004,Hairer2006Geometric} and references therein). The resulting variational time integrators are, for a nondissipative and nonforced case, symplectic and momentum preserving. Experience has shown that, while not being exactly energy conserving, constant time-step symplectic-momentum methods tend to exhibit good energy behavior over exponentially long times (see, for example, \cite{Kane2000NewmarkAlgorithm} for Lagrangian dynamics on $\mathrm{E}(3)^n$, \cite{Lee2007,Leok2011,Leok2012} for Lagrangian dynamics on $\mathrm{SE}(3)^n$, and references therein). 

It bears emphasis that variational time integrators for full body dynamics based on Lie group computations, such as the Lie group variational integrators developed by Lee, Leok and McClamroch \cite{Lee2007}, require the solution, in the space of rotations, of an implicit equation for the relative attitude between two integration steps. Alternatively, variational integrators for rigid body dynamics in terms of unit quaternions, such as the work by Wendlandt and Marsden \cite{Wendlandt1997}, require the enforcement constraints through Lagrange multipliers. Recently, Campello \cite{Campello2015} proposed an interesting parametrization of the space of rotations that uses exponential local coordinates represented by a rescaled form of the Rodrigues rotation vector. This parametrization eliminates the need of reprojections, constraints, and the use of complicated trigonometric or transcendental expressions. The time integration scheme proposed by Campello, however, is not derived from a discrete Lagrangian, it is implicit, in general, and its explicit version results in a first-order time integrator with stability issues at large time steps.

In this work, we combine and extend the work by Leok and Campello, and their co-workers, to develop an explicit, second-order, variational time integrator for full body dynamics that preserves the momenta of the continuous dynamics, such as linear and angular momenta, and exhibits near-conservation of total energy over exponentially long times.

The paper is organized as follows. The variational integrator based on rescaled Rodrigues parameters is defined in Section 2. We next illustrate the convergence, accuracy and conservation properties of the proposed time-stepping scheme with two selected examples of application. Fist, the dynamics of a three-dimensional pendulum in $\mathrm{SO}(3)$ is studied in Section 3. Second, the impact of a particle-binder torus against a rigid wall is studied in Section 4 using particle-binder-particle interactions in $\mathrm{SE}(3)^n$. Finally, a summary and concluding remarks are collected in Section 5.

\section{Variational integrator based on rescaled Rodrigues parameters}

The Lagrangian dynamics of a translating and rotating rigid body in three-dimensions, i.e., the Euclidean motion of a rigid body $\mathrm{SE}(3):(R,x)$, is given by $t \rightarrow (R, x, \dot{R},\dot{x})$, where $R \in \mathrm{SO}(3)$ is the attitude of the rigid body that transforms a vector from the body-fixed frame to the inertial frame, $x \in \mathbb{R}^3$ is the position vector of the origin of the body-fixed frame expressed in the inertial frame, $\dot{R}=R S(\omega)$ with $\omega$ being the angular velocity vector in the body-fixed frame, $S(\omega)$ is the skew matrix of vector $\omega \in \mathbb{R}^3$, and $\dot{x}$ is the translational velocity vector of the rigid body expressed in the inertial frame. We restrict attention to Lagrangian functions with a quadratic kinetic energy, that is defined by
\begin{equation}
L\left( R, x, \dot{R},\dot{x} \right)
=
\frac{1}{2} \dot{x}^T  m  \dot{x}
+
\frac{1}{2} \mathrm{tr}(R^T \dot{R} J_d \dot{R}^T R)
-
U(R,x)
\end{equation}
and, alternatively, for compactness, by
\begin{equation}
	\hat{L}\left( R, x, \omega, \dot{x}  \right)
	=
	L\left( R, x, R S(\omega), \dot{x} \right)	
	=
	\frac{1}{2} \dot{x}^T m \dot{x} 
	+
	\frac{1}{2} \omega^T J \omega
	-
	U(R,x)
\end{equation}
where $U(R,x)$ is the potential energy of the system, $m$ is the mass of the body, $J$ is the body's rotational moment of inertia in the body-fixed frame, with $J_d=\frac{1}{2}\mathrm{tr}(J) I - J$ and, thus, $J = \mathrm{tr}(J_d) I - J_d$. We also restrict attention to spherical bodies and the special case in which the origin of the body-fixed frame is located at the center of mass of the body. Therefore, the rotational moments of inertia are the scalars $J$ and $J_d=J/2$. Lagrangians of this form arise in a number of areas of application including particulate systems composed of spherical particles.

The configuration manifold $\mathrm{SE}(3)$ is a Lie group and its representation poses a special challenge in that the space of finite three-dimensional rotations is not linear. A common approach for this representation is to introduce a local linear approximation and local coordinates in the form of, among others, three angle coordinates (e..g., Euler angles) or exponential coordinates (e.g., Euler rotation vector $\theta \in \mathbb{R}^3$ such that $R=\mathrm{e}^{S(\theta)}$), which in turn introduce complexity in the analysis and computations. In this work,  we parametrize the space of rotations using exponential local coordinates and we choose a rescaled form of the Rodrigues rotation vector $\alpha \in \mathbb{R}^3$ recently proposed by Piment, Campello and Wriggers \cite{Pimenta2008, Campello2011}, and adopted by Campello and co-workers \cite{Campello2015} for the dynamics of granular systems, that reduces computational overhead. Henceforth, we express the rotation motion $R$ as a function of rotation vector $\alpha$. Specifically, the rescaled Rodrigues rotation vector is related to the Euler rotation vector by the following expression
\begin{equation}
\alpha = \frac{\tan(\|\theta\|/2)}{\|\theta\|/2} \theta
\end{equation}
and it is related to the rotation matrix $R(\alpha)$ by
\begin{equation}
R(\alpha)
=
I
+
\frac{4}{4+\| \alpha \|^2}
\left[  S(\alpha) + \frac{1}{2}  S(\alpha)^2 \right]	
\end{equation}
where the exponential map $R=\mathrm{e}^{S(\theta)}$ has been replaced by the closed form expression known as the Euler-Rodrigues rotation formula. This parametrization of the space of rotations is not global, in the sense that it is not able to represent rotations of magnitude $\| \theta \| = \pm \pi$. Similarly, the modified Rodrigues rotation vector \cite{Schaub2003}, i.e., $\tan(\|\theta\|/4)~\theta/\|\theta\|$, is singular at $\| \theta \| = \pm 2\pi$. On the contrary, the representation of the associated rotation matrix $R(\alpha)$ is free from trigonometric functions and the composition of two successive rotations $\alpha_{(1)}$ and $\alpha_{(2)}$ is determined by the following simple expression
\begin{equation}
\alpha_{(1)} \oplus \alpha_{(2)}
=
\frac{4}{4-\alpha_{(1)} \cdot \alpha_{(2)}}
\left[  
\alpha_{(1)}
+ 
\alpha_{(2)}
- 
\frac{1}{2} \alpha_{(1)} \times \alpha_{(2)}  \right]
\end{equation}
where  $\alpha_{(1)}$ and $\alpha_{(2)}$ are rescaled Rodrigues rotation vectors. It is worth noting that this composition automatically generates a rotation in $\mathrm{SO}(3)$, without the need of constraints or reprojection techniques.

\subsection{Discrete Lagrangian based on rescaled Rodrigues parameters}

Following the work of Lee, Leok and McClamroch \cite{Lee2007}, a variational time integrator is derived from the following discrete Lagrangian function
\begin{equation}
	\begin{aligned}		
		L_d(x_k,x_{k+1},\alpha_k,\alpha_k \oplus \Delta{\alpha}_{k}) 
		=
		\frac{h}{2}
		L\left(R(\alpha_k), x_k, \frac{R(\alpha_k \oplus \Delta\alpha_k)- R(\alpha_k)}{h},\frac{x_{k+1}-x_{k}}{h} \right)
		&+ \\
		\frac{h}{2}
		L\left(R(\alpha_k \oplus \Delta\alpha_k), x_{k+1}, \frac{R(\alpha_k \oplus \Delta\alpha_k)-R(\alpha_k)}{h},\frac{x_{k+1}-x_{k}}{h} \right)	&
	\end{aligned}
\end{equation}
and the following approximation of the action integral
\begin{equation}
	\mathcal{S} = \sum_{k=0}^{N-1} L_d(x_k,x_{k+1},\alpha_k,\alpha_k \oplus \Delta{\alpha}_{k}) 
	\approx
	\int_{t_0}^{t_N} L\left( R, x, \dot{R},\dot{x} \right) dt
\end{equation}
where $h$ is the time step, and $x_k$ and $\alpha_k$ are the position in the inertial frame and the rescaled Rodrigues rotation vector, respectively, at time $t_k=kh$. Variations are chosen to respect the geometry of the configuration space $SE(3):(R,x)$. Specifically, the variation of $x_k$ is given by $x_k^\epsilon=x_k+\epsilon \delta x_k+\mathcal{O}(\epsilon^2)$, where $\delta x_k \in \mathbb{R}^3$ vanishes at $k=0$ and $k=N$. Similarly, the variation of $R_k$ is given by $R^\epsilon_k = R_k \mathrm{e}^{\epsilon \eta_k}=R_k+\epsilon \delta R_k +\mathcal{O}(\epsilon^2)$, and thus $\delta R_k = R_k \eta_k$, where $\eta_k\in \mathrm{so}(3)$ is a skew-symmetric matrix that vanishes at $k=0$ and $k=N$. The discrete version of Hamilton’s variational principle then gives
\begin{equation}
	\begin{aligned}		
		\delta\mathcal{S} = 0 = 
		\sum_{k=1}^{N-1}
		&-\delta x_k \left[ \frac{1}{h} m (x_{k+1} - 2 x_k + x_{k-1}) + h \frac{\partial U_k}{\partial x_k}  \right] +
		\\
		&\mathrm{tr}\left( \eta_k \left[ \frac{1}{h} R_k^T R_{k+1} J_d - \frac{1}{h} J_d R_{k-1}^T R_k + h R_k^T \frac{\partial U_k}{\partial R_k}  \right]  \right)
	\end{aligned}
\end{equation}
with $R_k = R(\alpha_k)$, $R_{k+1}=R(\alpha_k \oplus \Delta\alpha_k)$, and $U_k=U(R_k,x_k)$.
It is worth noting that $R(\alpha_k \oplus \Delta\alpha_k)$ enforces $R_{k+1} \in \mathrm{SO}(3)$ automatically, without the need of constraints or reprojection techniques---cf. \cite{Lee2007}. Therefore, the following discrete Euler-Lagrange equations are obtained
\begin{equation}
	\frac{1}{h} m (x_{k+1} - 2 x_k + x_{k-1}) + h \frac{\partial U_k}{\partial x_k} =0
\end{equation}
and, since $\eta_k$ is skew-symmetric, 
\begin{equation}
	\frac{1}{h} \left( R_k^T R_{k+1} J_d - J_d R_{k-1}^T R_k 
	- J_d R_{k+1}^T R_k + R_k^T R_{k-1} J_d	\right)
	+ h \left( R_k^T \frac{\partial U_k}{\partial R_k} - \frac{\partial U_k}{\partial R_k}^T R_k  \right) = 0
\end{equation}
for any choice of $\delta x_k$ with $\delta x_0 = \delta x_N =0$, and any choice of $\eta_k$ that vanishes at $k=0$ and $k=N$.

By defining the momentum at time step $k$ as the conjugate of $x_k$ using the discrete Legendre transformation \cite{Lee2007}, the following expression is obtained
\begin{equation}
\begin{aligned}
	m v_k
	=
	- \frac{\partial L_d(x_k,x_{k+1},\alpha_k,\alpha_k \oplus \Delta{\alpha}_{k})}{\partial x_k}
	=
	\frac{m}{h} (x_{k+1}-x_k)+\frac{h}{2} \frac{\partial U_k}{\partial x_k}
\end{aligned}
\end{equation}
and, by employing the first discrete Euler-Lagrange equation obtained above, the momentum at time step $k+1$ is given by
\begin{equation}
\begin{aligned}
	m v_{k+1}
	&=
	- \frac{\partial L_d(x_{k+1},x_{k+2},\alpha_{k+1},\alpha_k \oplus \Delta{\alpha}_{k+1})}{\partial x_{k+1}}
	=
	\frac{\partial L_d(x_k,x_{k+1},\alpha_k,\alpha_k \oplus \Delta{\alpha}_{k})}{\partial x_{k+1}}
	\\
	&=
	\frac{m}{h} (x_{k+1}-x_k)-\frac{h}{2} \frac{\partial U_{k+1}}{\partial x_{k+1}}
\end{aligned}
\end{equation}
The discrete update map $(x_k , v_k) \rightarrow (x_{k+1} ,  v_{k+1})$ is then given by
\begin{subequations}
	\begin{align}
	x_{k+1} &= x_k +h v_k + \frac{h^2}{2m} F_k
	\\
	v_{k+1} &= v_k + \frac{h}{2m} (F_k + F_{k+1})
	\end{align}
\end{subequations}
where $F_k=-\partial U_k/\partial x_k$. We note that for $U(R,x)$, both $\alpha_k$ and $\alpha_{k+1}$ are needed to evaluate $F_k$ and $F_{k+1}$ and, thus, this is presented next in turn.

Similarly, using the discrete analogue of the Legendre transformation \cite{Lee2007}, we define the angular velocity vector in the body-fixed frame $\omega_k$ at time step $k$ as
\begin{equation}
	\begin{aligned}
	S(J \omega_k)
	=&
	- \frac{\partial L_d(x_k,x_{k+1},\alpha_k,\alpha_k \oplus \Delta{\alpha}_{k})}{\partial R_k} + \frac{\partial L_d(x_k,x_{k+1},\alpha_k,\alpha_k \oplus \Delta{\alpha}_{k})}{\partial R_k}^T
	\\
	=&\frac{1}{h} \left[ R_k^T R_{k+1} J_d - J_d R_{k+1}^T R_k \right] 
	- \frac{h}{2} S(M_k)
	\end{aligned}
\end{equation}
where $S(M_k)=(\partial U_k/\partial R_k)^T R_k -R_k^T (\partial U_k/\partial R_k)$. Thus, the discrete map $(\alpha_k , \omega_k) \rightarrow (\Delta \alpha_{k})$ is implicitly given by the following quadratic equation
\begin{equation}
h \left(J \omega_k + \frac{h}{2} M_k \right)
=
4 J R_k^T \frac{\Delta\alpha_k}{4+\|\Delta\alpha_k\|^2}
\end{equation}
with solution give by
\begin{equation}
	\Delta\alpha_{k}
	=
	\frac{2 h}{1+\sqrt{1- h^2 \| \omega_k + \frac{h}{2J} M_k \|^2}}
	R_k 
	\left[ \omega_k + \frac{h}{2J} M_k \right]		
\end{equation}
for $h^2 \| \omega_k + h M_k/2J \|<1$, i.e., for incremental rotations smaller than $\pi/2$. It is worth noting that the second solution of the quadratic equation doesn't fulfill the conditions of $\| \Delta \alpha_k \| \rightarrow 0$ for $\| \omega_k + h M_k/2J \| \rightarrow 0$. The systematic investigation of the discrete map $(\alpha_k , \omega_k) \rightarrow (\Delta \alpha_{k})$ for a general moment of inertia $J\in\mathbb{R}^3$, i.e., for non-spherical particles, and the elucidation of closed-form analytical solutions, are worthwhile directions of future research.

By the same token, using the discrete analogue of the Legendre transform, we define $\omega_{k+1}$ as
\begin{equation}
	\begin{aligned}
	S(J \omega_{k+1})
	=&
	\frac{\partial L_d(x_k,x_{k+1},\alpha_k,\alpha_k \oplus \Delta{\alpha}_{k})}{\partial R_{k+1}} - \frac{\partial L_d(x_k,x_{k+1},\alpha_k,\alpha_k \oplus \Delta{\alpha}_{k})}{\partial R_{k+1}}^T
	\\
	=&\frac{1}{h} \left[ J_d R_k^T R_{k+1} - R_{k+1}^T R_k J_d \right] 
	+ \frac{h}{2} S(M_{k+1})
	\end{aligned}
\end{equation}
and, using the second discrete Euler-Lagrange equation obtained above, the angular velocity vector at time step $k+$ is then given by
\begin{equation}
	\omega_{k+1} 
	= 
	R_{k+1}^T R_k  \left( \omega_k + \frac{h}{2J} M_k \right) + \frac{h}{2J} M_{k+1}
\end{equation}

Lastly, since the potential energy depends solely on the position and orientation of the rigid body, i.e., $U(\alpha, x)$, the following discrete Lagrangian map is explicit
\begin{equation}
(\alpha_k, x_k , \omega_k , v_k) \rightarrow (\alpha_{k+1} , x_{k+1} , \omega_{k+1} , v_{k+1})
\end{equation}
and the attitude of the rigid body $R(\alpha_k)$ can be reconstructed at any given time step.

It is worth noting that the Lie group variational integrators developed by Lee, Leok and McClamroch \cite{Lee2007} require the solution, in the space of rotations, of an implicit equation for the relative attitude between two integration steps. Here, the relative attitude between two integration steps is conveniently given by $\Delta{\alpha}_{k}$, which follows from an explicit, closed-form expression. As noted above, the generalization of these variational integrators based on rescaled Rodrigues parameters to granular systems of non-spherical particles remains as an interesting topic of study.

\subsection{Explicit variational time integrator}

The variational time itegrator is then defined by the discrete Lagrangian map
\begin{equation}
(\alpha_k, x_k , \omega_k , v_k) \rightarrow (\alpha_{k+1} , x_{k+1} , \omega_{k+1} , v_{k+1})
\end{equation}
given by
\begin{subequations}
	\begin{align}	
		{x}_{k+1} 
		&=
		 {x}_k + h {v}_k + \frac{h^2}{2 m}  {F}_k
		\\
		\Delta{\alpha}_{k}
		&=
		\frac{2 h}{1+\sqrt{1- h^2 \|  {\omega}_k + \frac{h}{2J} {M}_k \|^2}}
		R_k
		\left[  {\omega}_k + \frac{h}{2J} {M}_k \right]	
		\text{\hspace{.25in}with\hspace{.1in}}
		{\alpha}_{k+1} 
		= 
		{\alpha}_k 
		\oplus 	
		\Delta{\alpha}_{k}	
		\\
		{v}_{k+1} 
		&= 
		{v}_k 
		+ \frac{h}{2 m} \left[{F}_k+{F}_{k+1}\right]		
		\\
		{\omega}_{k+1}
		&=
		R_{k+1}^T R_k   
		\left[ {\omega}_k +  \frac{h}{2J} {M}_k \right]
		+
		\frac{h}{2 J} {M}_{k+1}	 
	\end{align}
\end{subequations}
Herein lies one advantage of using the rescaled Rodrigues parameters. The discrete map is an algebraic update that does not require trigonometry, unlike expressions involving quaternions or Euler angles.

This explicit variational time integrator can be recast by adopting a inertially fixed frame and defining ${\Omega} = R {\omega}$ as the angular velocity vector in the inertial frame (as opposed to ${\omega}$ in the body-fixed frame), sometimes referred to as the spatial angular velocity vector. Using this change of variables, the variational time integrator is then defined by the following discrete Lagrangian map
\begin{equation}
(\alpha_k, x_k , \Omega_k , v_k) \rightarrow (\alpha_{k+1} , x_{k+1} , \Omega_{k+1} , v_{k+1})
\end{equation}
and given by
\begin{subequations}
	\label{Eqn-VTI-1}
	\begin{align}	
	{x}_{k+1} 
	&=
	{x}_k + h {v}_k + \frac{h^2}{2 m}  {F}_k
	\\
	\Delta{\alpha}_{k}
	&=
	\frac{2 h}{1+\sqrt{1- h^2 \| {\Omega}_k + \frac{h}{2J} R{M}_k \|^2}}
	\left[  {\Omega}_k + \frac{h}{2J} R{M}_k \right]		
	\text{\hspace{.25in}with\hspace{.1in}}
	{\alpha}_{k+1} 
	= 
	{\alpha}_k 
	\oplus 	
	\Delta{\alpha}_{k}	
	\\
	{v}_{k+1} 
	&= 
	{v}_k 
	+ \frac{h}{2 m} \left[{F}_k+{F}_{k+1}\right]		
	\\
	{\Omega}_{k+1}
	&=
	{\Omega}_k +  \frac{h}{2J} \left[ R{M}_k + R{M}_{k+1}	\right] 
	\end{align}
\end{subequations}
where $S(RM_k) = R_k S(M_k) R_k^T = R_k (\partial U_k/\partial R_k)^T - (\partial U_k/\partial R_k) R_k^T$. For completeness, we note that $\dot{R}_k=S(\Omega_k) R_k$.

\subsection{Second-order time integrator \textit{\`{a} la} explicit Newmark}

It is worth noting that the first two leading order terms in the expansion of the rescaled Rodrigues parameters update are given 
\begin{equation}
\Delta{\alpha}_{k}
=
\frac{2 h}{1+\sqrt{1- h^2 \| {\Omega}_k + \frac{h}{2J} R{M}_k \|^2}}
\left[  {\Omega}_k + \frac{h}{2J} R{M}_k \right]		
=
h{\Omega}_k + \frac{h^2}{2J} R{M}_k + \mathcal{O}\left( h^3 \right)	
\end{equation}
Therefore, the discrete Lagrangian map is simplified, by keeping linear and second order terms in $h$, as follows
\begin{subequations}
	\label{Eqn-VTI-2}	
	\begin{align}	
	{x}_{k+1} 
	&=
	{x}_k + h {v}_k + \frac{h^2}{2 m}  {F}_k
	\\
	\Delta{\alpha}_{k}
	&=
	h{\Omega}_k + \frac{h^2}{2J} R{M}_k	
	\text{\hspace{.5in}with\hspace{.1in}}
	{\alpha}_{k+1} 
	= 
	{\alpha}_k 
	\oplus 	
	\Delta{\alpha}_{k}
	\\
	{v}_{k+1} 
	&= 
	{v}_k 
	+ \frac{h}{2 m} \left[{F}_k+{F}_{k+1}\right]		
	\\
	{\Omega}_{k+1}
	&=
	{\Omega}_k +  \frac{h}{2J} \left[ R{M}_k + R{M}_{k+1}	\right] 
	\end{align}
\end{subequations}
The first and third equations in this map constitute the explicit Newmark time integrator ($\beta=0$ and $\gamma=1/2$), or the velocity Verlet algorithm, and the second and fourth equations exhibit the same mathematical structure. Therefore, we coin this variational integrator as a second-order time integrator \textit{\`{a} la} explicit Newmark or velocity Verlet.

\subsection{Other time integrators based on rescaled Rodrigues parameters}

The time integration scheme proposed by Campello and co-workers \cite{Campello2015} utilizes rescaled Rodrigues parameters and it is given by the following implicit discrete map
\begin{subequations}	
\begin{align}	
	{x}_{k+1} 
	&=
	{x}_k + h {v}_{k+1}
	\\
	\Delta{\alpha}_{k}
	&=
	h{\Omega}_{k+1}	
	\text{\hspace{.5in}with\hspace{.1in}}
	{\alpha}_{k+1} 
	= 
	{\alpha}_k 
	\oplus 	
	\Delta{\alpha}_{k}	
	\\
	{v}_{k+1} 
	&= 
	{v}_{k} + 
	\frac{h}{m} \left[\phi {F}_k+(1-\phi) {F}_{k+1}\right]		
	\\
	{\Omega}_{k+1}
	&=
	{\Omega}_{k}
	+
	\frac{h}{J} \left[ \phi R{M}_k + (1-\phi) R{M}_{k+1}  \right] 
\end{align}
\end{subequations}
For $\phi=1$, the integration scheme is explicity and given by
\begin{subequations}
\label{Eqn-VTI-3}	
\begin{align}	
	{x}_{k+1} 
	&=
	{x}_k + h {v}_{k+1}
	\\
	\Delta{\alpha}_{k}
	&=
	h{\Omega}_{k+1}	
	\text{\hspace{.5in}with\hspace{.1in}}
	{\alpha}_{k+1} 
	= 
	{\alpha}_k 
	\oplus 	
	\Delta{\alpha}_{k}
	\\
	{v}_{k+1} 
	&= 
	{v}_{k} + 
	\frac{h}{m} {F}_k
	\\
	{\Omega}_{k+1}
	&=
	{\Omega}_{k}
	+
	\frac{h}{J} R{M}_k
\end{align}
\end{subequations}
It is worth noting that this algorithm posses the structure of the symplectic Euler integrator and, thus, it is a first-order time integrator. Furthermore, it can be readily shown that it is a variational time integrator, as it is derived from the following discrete Lagrangian
\begin{equation}
	L^*_d(x_k,x_{k+1},\alpha_k,\alpha_k \oplus \Delta{\alpha}_{k}) 
	=
	h
	L\left(R(\alpha_k), x_k, \frac{R(\alpha_k \oplus \Delta\alpha_k)- R(\alpha_k)}{h},\frac{x_{k+1}-x_{k}}{h} \right)
\end{equation}
as follows
\begin{subequations}
\begin{align}
	m v_k
	&=
	- \frac{\partial L^*_d(x_k,x_{k+1},\alpha_k,\alpha_k \oplus \Delta{\alpha}_{k})}{\partial x_k}
	=
	\frac{m}{h} (x_{k+1}-x_k)+h \frac{\partial U_k}{\partial x_k}
	\\
	m v_{k+1}
	&=
	- \frac{\partial L^*_d(x_{k+1},x_{k+2},\alpha_{k+1},\alpha_k \oplus \Delta{\alpha}_{k+1})}{\partial x_{k+1}}
	=
	\frac{m}{h} (x_{k+1}-x_k)
	\\
	S(J \omega_k)
	&=
	- \frac{\partial L^*_d(x_k,x_{k+1},\alpha_k,\alpha_k \oplus \Delta{\alpha}_{k})}{\partial R_k} 
	+ \frac{\partial L^*_d(x_k,x_{k+1},\alpha_k,\alpha_k \oplus \Delta{\alpha}_{k})}{\partial R_k}^T
	\\
	&=\frac{1}{h} \left[ R_k^T R_{k+1} J_d - J_d R_{k+1}^T R_k \right] 
	- h S(M_k) \notag
	\\
	S(J \omega_{k+1})
	&=
	\frac{\partial L^*_d(x_k,x_{k+1},\alpha_k,\alpha_k \oplus \Delta{\alpha}_{k})}{\partial R_{k+1}} - \frac{\partial L^*_d(x_k,x_{k+1},\alpha_k,\alpha_k \oplus \Delta{\alpha}_{k})}{\partial R_{k+1}}^T
	\\
	&=\frac{1}{h} \left[ J_d R_k^T R_{k+1} - R_{k+1}^T R_k J_d \right] \notag
\end{align}
\end{subequations}
whence the discrete Lagrangian map proposed by Campello is recovered after retaining the leading order terms in $\Delta\alpha_{k}$, that is
\begin{equation}
	\Delta\alpha_{k}
	=
	\frac{2 h}{1+\sqrt{1- h^2 \| \Omega_{k+1} \|^2}}
	\Omega_{k+1} 
	=
	h{\Omega}_{k+1} +  \mathcal{O}\left( h^3 \right)	
\end{equation}

\subsection{Conservation properties and convergence}

The classical theorem of Noether states that if a Lagrangian is invariant under the lifted action of a Lie group, the associated momentum map is a constant of the motion, that is it is preserved along the trajectories. The body of work on variational time integrators (see \cite{Hairer2006Geometric} and references therein) indicates that if a discrete Lagrangian is invariant under the diagonal action of a symmetry group, a discrete Noether’s theorem holds, and the discrete flow preserves the discrete momentum map.  In addition, variational time integrators preserve the symplecticity of the discrete flow. Here we restrict attention to the Lie group $\mathrm{SE}(3)^n$, that is to the Euclidean motion of a system of $n$ rigid bodies.

{\bf Conservation of translational momentum.} If the discrete Lagrangian $L_d$ is invariant under translations $u \in \mathrm{E}(3)^n$, 
\begin{equation}		
	L_d(x_k+u,x_{k+1}+u,R_k,R_{k+1}) 
	=
	L_d(x_k,x_{k+1},R_k,R_{k+1}) 
\end{equation}
which holds true for $U(R,x+u)=U(R,x)$, then the discrete total translational momentum in the body-fixed frame is a constant of the discrete motion, that is
\begin{equation}
m_{(1)} v_{(1)k}+ ... + m_{(n)} v_{(n)k} = m_{(1)} v_{(1)0} + ... + m_{(n)} v_{(n)0}
\end{equation}

{\bf Conservation of angular momentum.}  If the discrete Lagrangian $L_d$ is invariant under rotations about the center of mass $\hat{R} \in \mathrm{SO}(3)^n$,
\begin{equation}
	L_d(x_k, x_{k+1}, \hat{R} R_k, \hat{R} R_{k+1}) 
	=
	L_d(x_k,x_{k+1},R_k,R_{k+1}) 
\end{equation}
which holds true for $U(\hat{R} R, x)=U(R,x)$, then the total angular momentum in the body-fixed frame is a constant of the discrete motion, that is
\begin{equation}
J_{(1)} \Omega_{(1)k} + ... + J_{(n)} \Omega_{(n)k}
=
J_{(1)} \Omega_{(1)0} + ... + J_{(n)} \Omega_{(n)0}
\end{equation}

{\bf Conservation of total energy.} While variational integrators do not exactly preserve energy, backward error analysis \cite{Hairer2006Geometric} shows that they preserve a modified Hamiltonian that is close to the original Hamiltonian for exponentially long times. In practice, the energy error is bounded and does not exhibit a drift.

{\bf Symmetry or time-reversibility.} Symmetry or time-reversibility follows directly from the definition of the discrete Lagrangian map and the adoption of rescaled Rodrigues vectors to represent rotations, which generates exact composition of rotations in $\mathrm{SO}(3)$.

{\bf Convergence.} The space of the discrete trajectories is the Sobolev space $H^1(t_0,t_N) := W^{1,2}(t_0,t_N)$ and, hence,  distances in the space of trajectories are measured with respect to the $H^1$-norm
\begin{equation}
\| q \|_{H^1(t_0,t_N)}^2
=
\int_{t_0}^{t_N}
\left(
\| x(t) \|^2
+
\| v(t) \|^2
+
\|\log(R(t))\|_F^2 /2
+
\|\Omega(t)\|^2
\right)
dt
\end{equation}
with the norm of the rotation $R$ defined as $\| \log(R) \|_F /\sqrt{2} = |\theta| = 2 \arctan(\|\alpha(t)\|/2)$ (see \cite{Huynh2009Metrics}). Hence, the error in the discrete total energy $E_h(t)$, for a time step of $h$, is measured with respect to the $H^0$-norm as follows
\begin{equation}
\| E_h(t) - E(0) \|_{H^0(t_0,t_N)}^2 = \int_{t_0}^{t_N} (E_h(t)-E_0)^2 dt 
\end{equation}
where $E_0$ is the initial kinetic and mechanical energy input to the system. For future reference, the relative error in trajectories and total energy are defined as follows
\begin{equation}
q(t)\text{-error}
=
\frac{\left| \| q_h \|_{H^1(t_0,t_N)}^2 - \| q \|_{H^1(t_0,t_N)}^2 \right|^{1/2}}{ \| q \|_{H^1(t_0,t_N)} }
\end{equation}
\begin{equation}
E(t)\text{-error} = \frac{\| E_h(t) - E(0) \|_{H^0(t_0,t_N)}}{\| E(0) \|_{H^0(t_0,t_N)}} 
\end{equation}

\section{Three-dimensional pendulum in $\mathrm{SO}(3)$}

A three-dimensional pendulum is a rigid body of mass $m$ and rotational moment of inertia $J$ supported by a fixed pivot and acted on by gravity $g$. An inertial Euclidean frame is selected so that the first two axes lie in a horizontal plane and the third axis is vertical, and thus gravity is $-g~e_3$. The origin of the inertial Euclidean frame is selected to be the location of the pendulum pivot $\rho_0$. The potential energy is given by
\begin{equation}
U(R) = - m g~e_3 \cdot R \rho_0
\end{equation}
with $\rho_0=[0,0,1]$ and $e_3=[0,0,1]$. The moment is then
\begin{equation}
RM
=R \rho_0 \times R^T e_3
\end{equation}
and the conserved quantities are
\begin{subequations}
\begin{align}	
E = \frac{1}{2} J \| \Omega \|^2 -m g~e_3 \cdot R \rho_0 =&~\mathrm{const}
\\
R \rho_0 \cdot \Omega =&~\mathrm{const}
\\
\| R \rho_0 \| =&~\mathrm{const}
\end{align}
\end{subequations}

\begin{figure}[htbp]
	\centering{
		\begin{tabular}{cc}
			\includegraphics[clip, trim=0cm 6cm 0cm 6cm, scale=0.32]{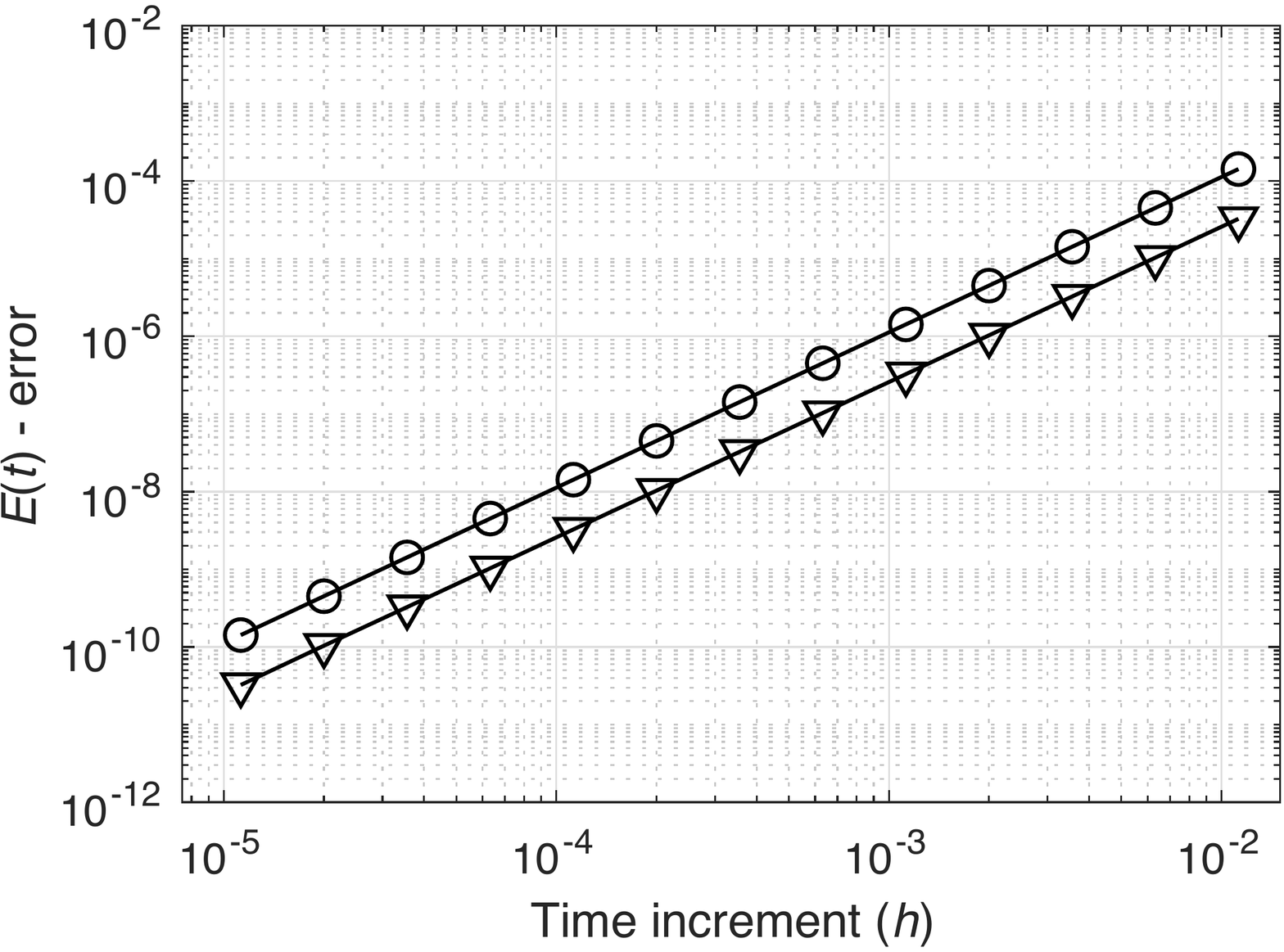}
			&
			\includegraphics[clip, trim=0cm 6cm 0cm 6cm, scale=0.32]{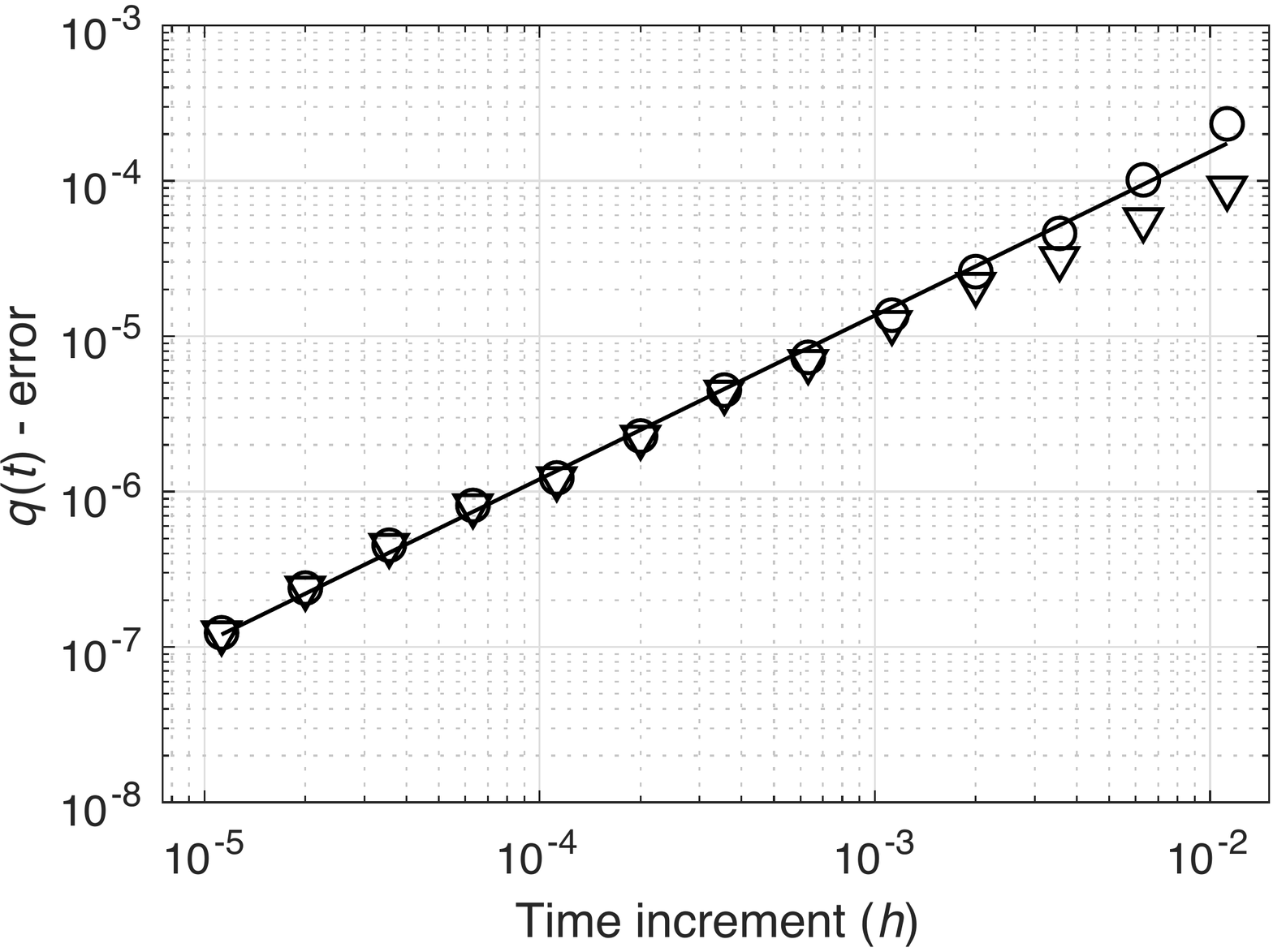}
			\\
			\includegraphics[clip, trim=0cm 6cm 0cm 6cm, scale=0.32]{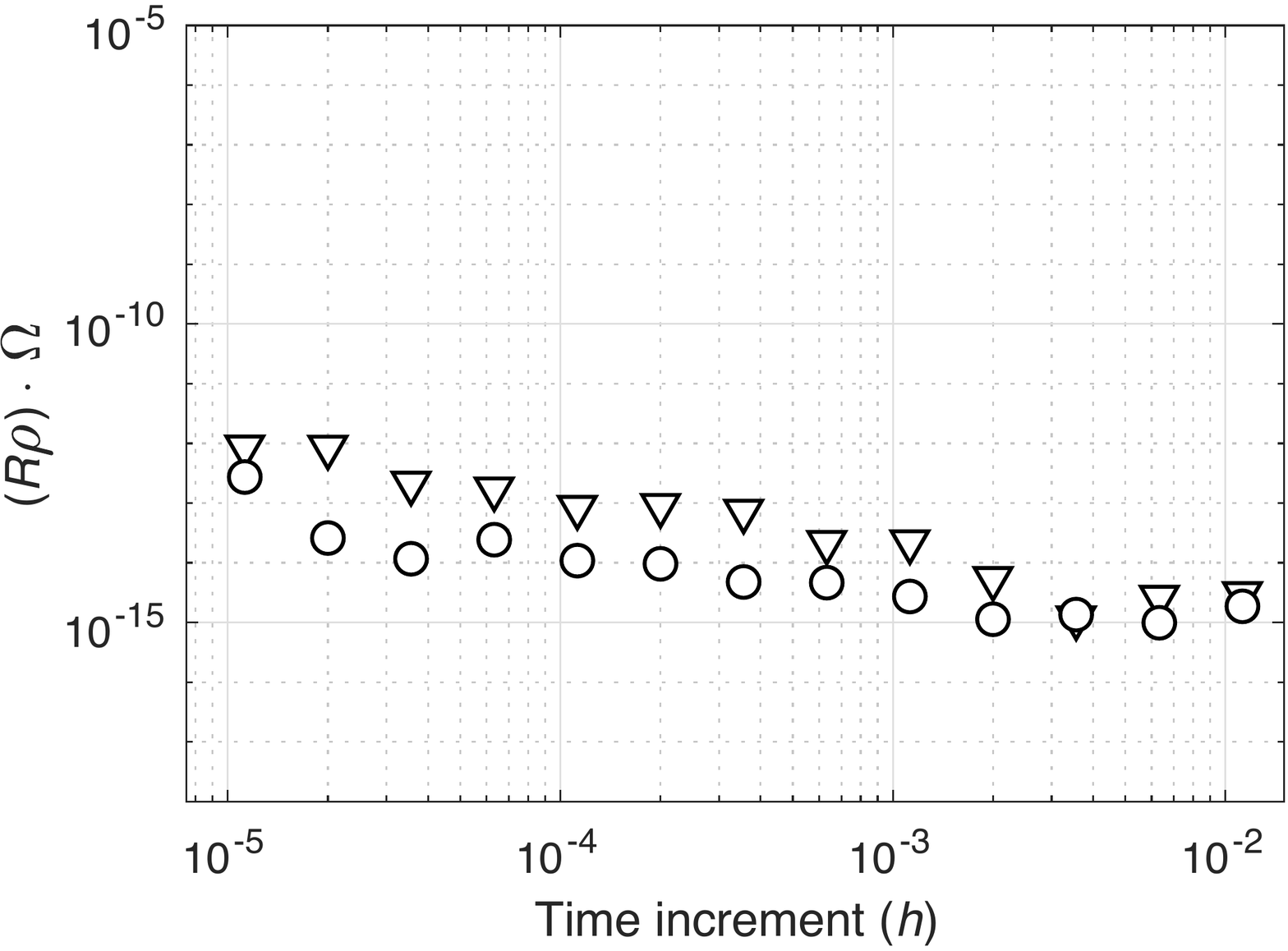}	
			&
			\includegraphics[clip, trim=0cm 6cm 0cm 6cm, scale=0.32]{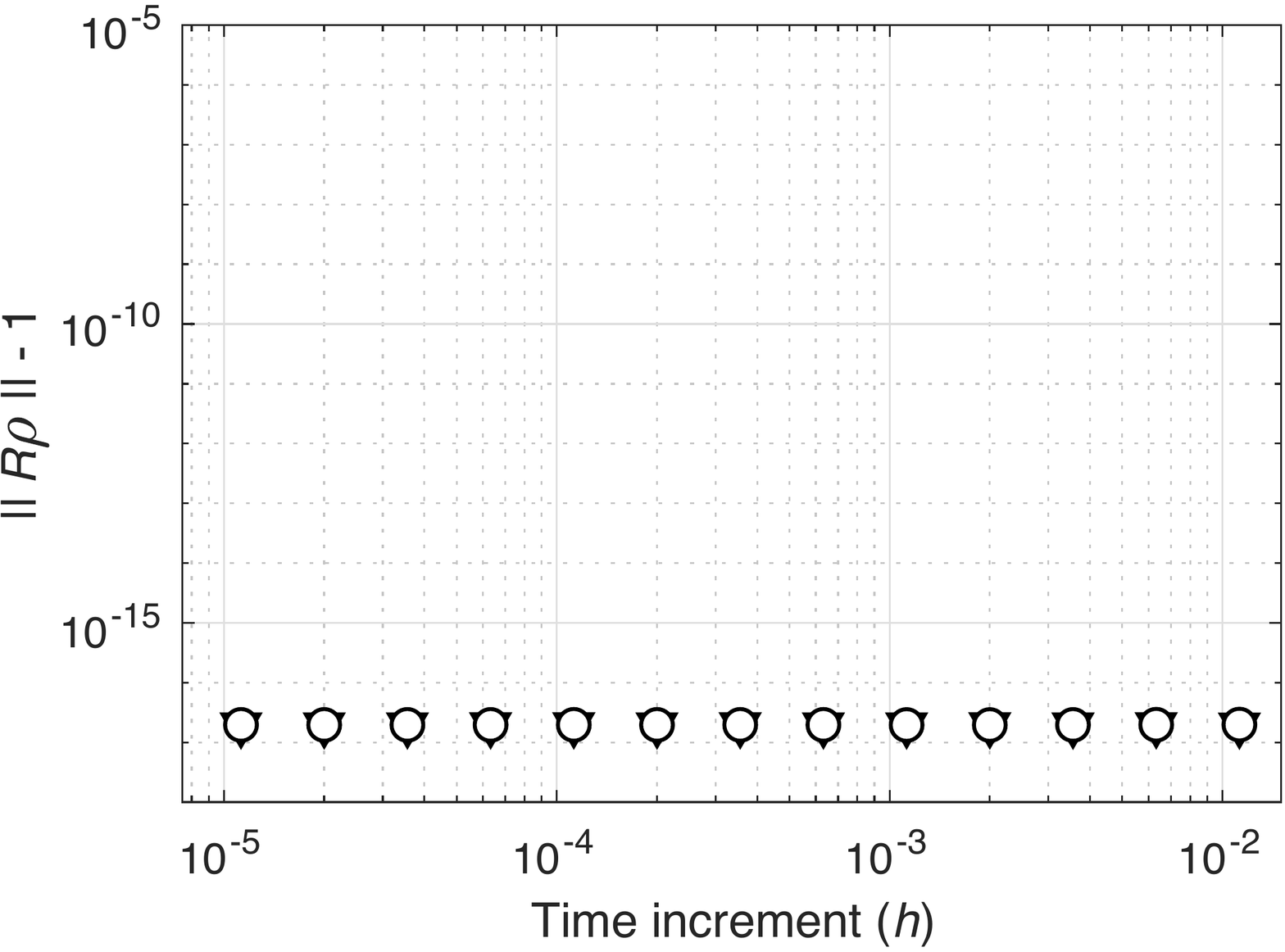}	
		\end{tabular}
	}
	\caption{Both explicit time integrators---discrete Lagrangian map \eqref{Eqn-VTI-1} in $\Circle$, and discrete Lagrangian map \eqref{Eqn-VTI-2} in $\triangledown$---exhibit (i) quadratic convergence of energy in the $H^0$-norm using the $E(t)\text{-error}$, (ii) linear convergence of trajectories in the $H^1$-norm using the $q(t)\text{-error}$, (iii) exact conservation, to machine precision, of $R \rho_0 \cdot \Omega=0$, and (iv) exact conservation, to machine precision, of $\| R \rho_0 \| = 1$. 
	}
	\label{Fig-3DPendulum}
\end{figure}

We adopt the following initial conditions to showcase the conservation, accuracy, longterm behavior and convergence properties of the second-order, explicit time integrator based on rescaled Rodrigues parameters,
\begin{equation}
\begin{aligned}
{\alpha}(t=0) =&~[0,1,0]~2~\mathrm{tan}(3\pi/8)
\\
{\Omega}(t=0) =&~[1,0,1]~0.4~\mathrm{sin}(\pi/4)^2
\end{aligned}
\end{equation}
Figure~\ref{Fig-3DPendulum} shows (i) quadratic convergence of energy in the $H^0$-norm using the $E(t)\text{-error}$, (ii) linear convergence of trajectories in the $H^1$-norm using the $q(t)\text{-error}$, (iii) exact conservation, to machine precision, of $R \rho_0 \cdot \Omega=0$, and (iv) exact conservation, to machine precision, of $\| R \rho_0 \| = 1$. It is interesting to note that the two explicit time integrators developed, that is the discrete Lagrangian maps given by \eqref{Eqn-VTI-1} and \eqref{Eqn-VTI-2}, show the same conservation and convergence properties.

\section{Particle-binder-particle interactions in $\mathrm{SE}(3)^n$}

\begin{figure}[H]
	\centering{
		\includegraphics[scale=.30]{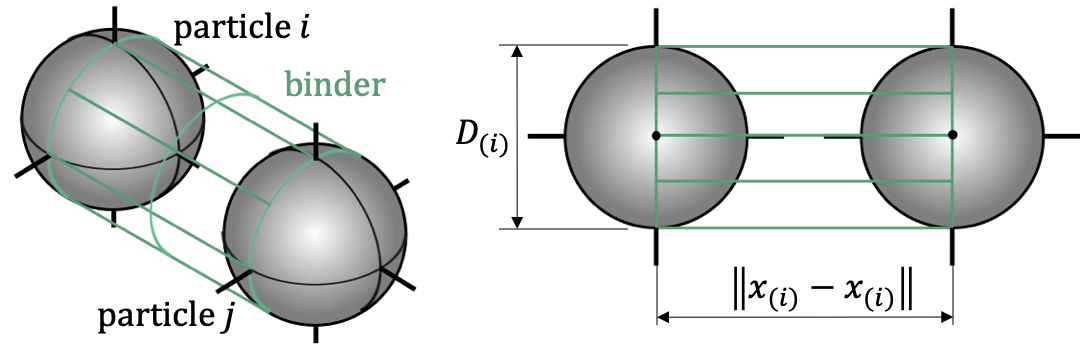}
	}
	\caption{Particles-binder-particles interaction.}
	\label{Fig-Particle-Binder-Particle}
\end{figure}

Particle-binder-particle interactions are used in particle-scale models of particle-binder composites, i.e, of composite materials that consist of a large concentration of hard particles, called fillers,
randomly dispersed in the matrix of a soft material. These composites exhibit strong nonlinearities in their static and dynamic mechanical response, owing to the presence of phases with largely different stiffness and complex relative motion between hard particles connected by soft binder. Specifically, the total strain potential energy in the system can be described by particles-binder-particles interactions (see Fig. \ref{Fig-Particle-Binder-Particle}) that account for binder torsional and bending loads ($U^\mathrm{m}$), binder axial loads ($U^\mathrm{a}$), binder shear loads ($U^\mathrm{s}$) and, when binder is ruptured under compression, particle-particle axial loads ($U^\mathrm{pp}$). The total strain potential energy is then given by 
\begin{equation}
U(R,x;R_0,x_0) 
= 
U^\mathrm{m}(R;R_0) 
+ 
U^\mathrm{a}(x;x_0) 
+ 
U^\mathrm{s}(R,x;R_0,x_0) 
+ 
U^\mathrm{pp}(x;x_0)
\end{equation}
where $x_0=x(t=0)$ and $R_0=R(t=0)$. Particle-wall interaction are described as a special case of particle-particle interactions.

The study of the dynamic response of these particle-binder composites, such as polymer-bonded explosives \cite{agarwal2020effects,agarwal2021constitutive,Baker2021,Poorsolhjouy2021} and asphaltic concrete \cite{Zhu1996, Si2002}, requires time integrators with good long-term behavior, good conservation properties, and high computational efficiency in the treatment of rotational degrees of freedom. We therefore use the impact of a particle-binder torus on a rigid wall as an example to illustrate the performance of the two second-order explicit variational time integrators developed in this work, i.e., the discrete Lagrangian maps given by \eqref{Eqn-VTI-1} and \eqref{Eqn-VTI-2}, and the first-order variational time integrator given by the discrete Lagrangian map \eqref{Eqn-VTI-3}. 

Next, we show expressions for each strain energy and their corresponding force and moment, in terms of the rescaled Rodrigues rotation vectors. For simplicity, we assume that $R_0=I$ and that the stiffness coefficients of the binder ($K^\mathrm{m}_{(ij)}$, $K^\mathrm{a}_{(ij)}$, $K^\mathrm{s}_{(ij)}$) and of the particle ($K^\mathrm{pp}_{(ij)}$) are constants.

{\bf Binder torsional and bending loads.} The strain energy associated to binder torsional and bending loads, i.e., $U^\mathrm{m}$, is given by
\begin{equation}	
	U^\mathrm{m} 
	=
	\frac{1}{2} \sum K^\mathrm{m}_{(ij)} \| \theta_{(ij)} \|^2
	=
	2 \sum K^\mathrm{m}_{(ij)} \arctan(\|\alpha_{(ij)}\|/2)^2
\end{equation}
where $\theta_{(ij)}$ is the relative Euler angle, and  $\alpha_{(ij)}$ is the relative rescaled Rodrigues vector, between particles $(i)$ and $(j)$, that is 
\begin{equation}
R_{(ij)}=\mathrm{e}^{S(\theta_{(ij)})}=R_{(i)} R_{(j)}^T
\hspace{.3in}\text{and}\hspace{.3in}
\alpha_{(ij)} = -\alpha_{(i)} \oplus \alpha_{(j)}
\end{equation}
with  $R_{(i)}= R(\alpha_{(i)})$ and $R_{(j)}=R(\alpha_{(j)})$. By definition, the moment is then given by
\begin{equation}
S(RM^\mathrm{m}_{(i)}) 
= 
R_{(i)} S(M^\mathrm{m}_{(i)}) R_{(i)}^T 
= 
R_{(i)} 
\left(\frac{\partial U^\mathrm{m}}{\partial R_{(i)}}\right)^T 
- 
\frac{\partial U^\mathrm{m}}{\partial R_{(i)}}
R_{(i)}^T 
\end{equation}
and, for a system of two particles $(1)$ and $(2)$, it simplifies to
\begin{equation}
RM_{(1)}
=
-RM_{(2)}
=
- K^\mathrm{m}_{(12)}  \theta_{(12)}
=
K^\mathrm{m}_{(12)} 
2 \arctan(\|\alpha_{(12)}\|/2) \frac{\alpha_{(12)}}{\|\alpha_{(12)}\|}
\end{equation}
which can be efficiently computed using the rescaled Rodrigues vectors---cf. \cite{Wang2009}. The derivation of this closed form expression follows from
\begin{equation}
	\begin{aligned}	
	U^\mathrm{m} 
	=& 
	\frac{1}{2} K^\mathrm{m}_{(12)} \| \theta_{(12)} \|^2
	=
	\frac{1}{2} K^\mathrm{m}_{(12)} 
	\left[
	\arccos\left(\frac{\mathrm{tr}(R_{(12)})-1}{2}\right)
	\right]^2
	\\	
	\frac{\partial U^\mathrm{m}}{\partial R_{(i)}}
	=&
	-
	\frac{1}{2} K^\mathrm{m}_{(12)}
	\| \theta_{(12)} \|
	\left[ 1 - \left(\frac{\mathrm{tr}(R_{(12)})-1}{2}\right)^2  \right]^{-1/2}
	\frac{\partial \mathrm{tr}(R_{(12)})}{\partial R_{(i)}}
	\end{aligned}
\end{equation}
with
\begin{equation}
\frac{\partial \mathrm{tr}(R_{(12)})}{\partial R_{(1)}}
=
R_{(2)}
\hspace{.3in}\text{and}\hspace{.3in}
\frac{\partial \mathrm{tr}(R_{(12)})}{\partial R_{(2)}}
=
- R_{(2)}  R_{(1)}^T  R_{(2)}
\end{equation}
and, hence, from
\begin{equation}
S(RM_{(1)}) 
=
-S(RM_{(2)}) 
=
-
\frac{K^\mathrm{m}_{(12)} \|\theta_{(12)}\|}{2 \sin(\|\theta_{(12)}\|)}
\left(
R_{(12)}-R_{(12)}^T
\right)
=
- K^\mathrm{m}_{(12)}
S(\theta_{(12)})
\end{equation}

{\bf Binder axial loads.} The strain energy associated to binder axial load, i.e., $U^\mathrm{a}$, is given by
\begin{equation}
	U^\mathrm{a} = \frac{1}{2} \sum K^\mathrm{a}_{(ij)} \left[ \frac{\| x_{(i)} - x_{(j)}\| }{\| x_{(i)0} - x_{(j)0}\|} - 1  \right]^2 
\end{equation}
and, for a system of two particles, the force is given by
\begin{equation}
	F_{(1)} 
	=
	-\frac{\partial U^\mathrm{a}}{\partial x_{(1)}} =-\frac{K^{a}_{(12)}}{\| x_{(1)0}-x_{(2)0} \|} \left[ \frac{\|x_{(1)} - x_{(2)} \|}{\| x_{(1)0}-x_{(2)0} \|} - 1 \right] 
	\frac{x_{(1)} - x_{(2)}}{\| x_{(1)} - x_{(2)} \|}
	= 
	-F_{(2)}
\end{equation}

{\bf Binder shear loads.} The strain energy associated to binder shear load, i.e., $U^\mathrm{s}$, is given by
\begin{equation}
U^\mathrm{s} 
= 
\frac{1}{2}
\sum K^\mathrm{s}_{(ij)} 
\left[ 
1 - \mathrm{tr}\left(\frac{R_{(i)}(x_{(i)0} - x_{(j)0})}{\| x_{(i)0}-x_{(j)0} \|} \otimes \frac{(x_{(i)} - x_{(j)}) R_{(j)}^T}{\| x_{(i)}-x_{(j)} \|}\right)
\right]^2
\end{equation}
and, for a system of two particles, the force is given by
\begin{equation}
	\begin{aligned}		
	F_{(1)} 
	=
	-\frac{\partial U^\mathrm{s}}{\partial x_{(1)}} 
	=	
	K^\mathrm{s}_{(12)}
	\left[ 
	1 - \mathrm{tr}\left(\frac{R_{(1)}(x_{(1)0} - x_{(2)0})}{\| x_{(1)0}-x_{(2)0} \|} \otimes \frac{(x_{(1)} - x_{(j)}) R_{(2)}^T}{\| x_{(1)}-x_{(2)} \|}\right)
	\right] 
	&
	\\
	\frac{ R_{(1)}(x_{(1)0} - x_{(2)0})R_{(2)} }
	       {\| x_{(1)} - x_{(2)} \| \| x_{(1)0} - x_{(2)0} \|}
	\left[ I
	-
	\frac{(x_{(1)} - x_{(2)}) \otimes (x_{(1)} - x_{(2)})}{\| x_{(1)} - x_{(2)} \|^2}	
	   \right]
	&= - F_{(2)}
	\end{aligned}
\end{equation}
and the moment by
\begin{equation}
	\begin{aligned}
	RM_{(1)} 
	= 
	K^\mathrm{s}_{(12)}
	\left[ 
	1 - \mathrm{tr}\left(\frac{R_{(1)}(x_{(1)0} - x_{(2)0})}{\| x_{(1)0}-x_{(2)0} \|} \otimes \frac{(x_{(1)} - x_{(2)}) R_{(2)}^T}{\| x_{(1)}-x_{(2)} \|}\right)
	\right]
	&
	\\
	\frac{ \left[ R_{(1)}(x_{(1)0}-x_{(2)0})\right] \times \left[(x_{(1)}-x_{(2)})  R_{(2)}^T \right]}{\| x_{(1)0} - x_{(2)0} \| \| x_{(1)}- x_{(2)} \|}
	&=
	- RM_{(2)}
	\end{aligned}
\end{equation}
Recasting these expressions in terms of rescaled Rodrigues vectors is desirable, if beyond the scope of this work.

{\bf Particle-particle axial loads.} The strain energy associated to particle-particle axial load, i.e., $U^\mathrm{pp}$, is given by
\begin{equation}
	U^\mathrm{pp} 
	= 
	\frac{5}{2} \sum K^\mathrm{pp}_{(ij)} 
	\left[1 - \frac{\| x_{(i)} - x_{(j)}\|}{D_{(i)}/2+D_{(j)}/2} \right]_{+}^{5/2} 
\end{equation}
with $[ \cdot ]_+=\max\{\cdot,0\}$ and, for a system of two particles, the axial force is given by
\begin{equation}
	F_{(1)} 
	=
	-\frac{\partial U^\mathrm{pp}}{\partial x_{(1)}} 
	=
	\frac{K^{pp}_{(12)}}{D_{(1)}/2+D_{(2)}/2} 
	\left[1 - \frac{\| x_{(1)} - x_{(2)}\|}{D_{(1)}/2+D_{(2)}/2} \right]_{+}^{3/2} 
	\frac{x_{(1)} - x_{(2)}}{\| x_{(1)} - x_{(2)} \|}
= 
-F_{(2)}
\end{equation}

\begin{figure}[H]
	\centering{
		\includegraphics[scale=0.65]{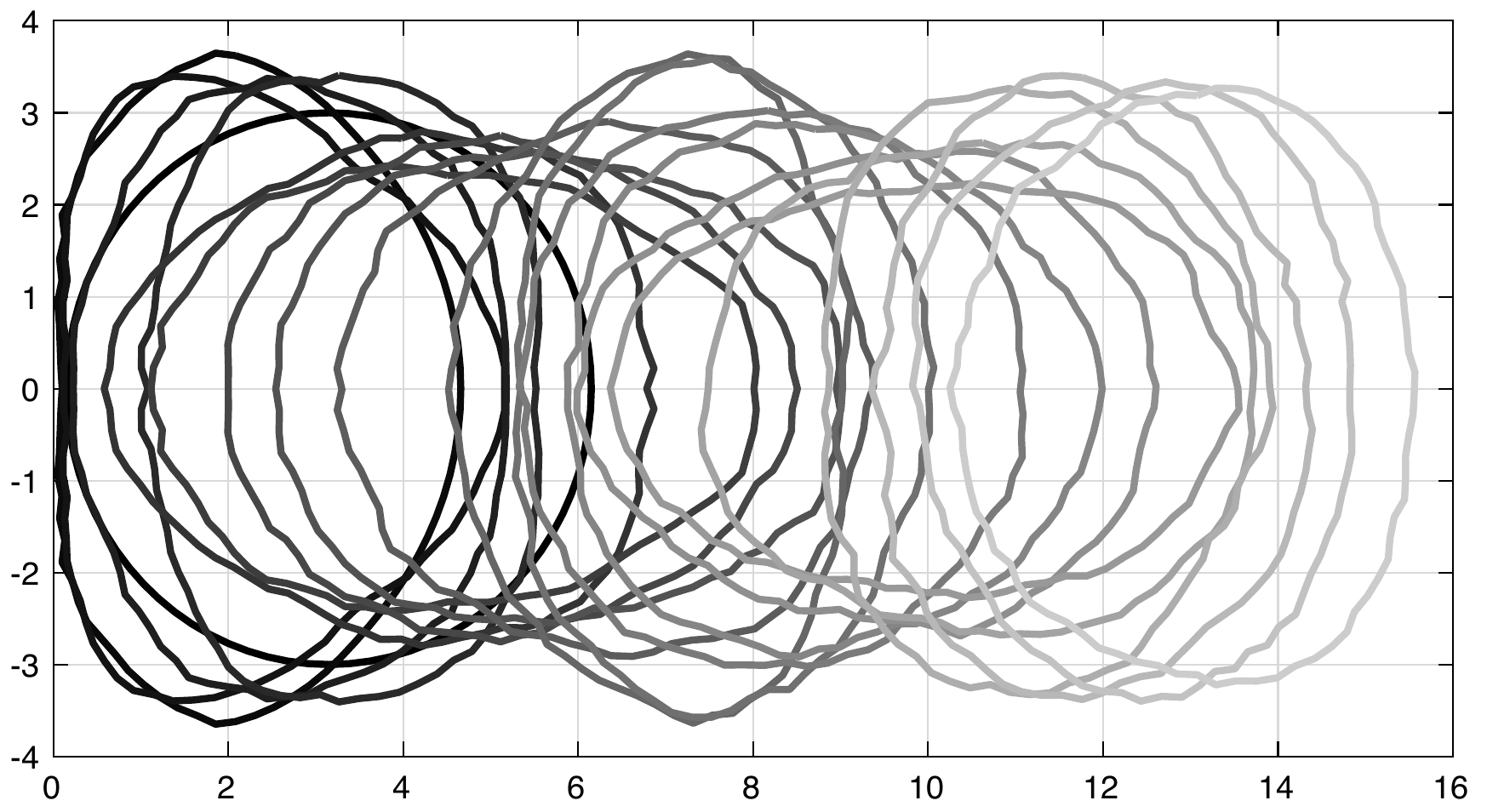}
	}
	\caption{Sequence of deformed states of the solid torus experiencing large deformation and vibrations after impact against a rigid wall with normal $[1,0,0]$ located at $x_w=0$. These numerical results are obtained using the second-order, explicit time integrator given by discrete map \eqref{Eqn-VTI-2}, with $h=10^{-3}$, $t_N=35$.}
	\label{Fig-TorusImpact-XY}
\end{figure}

\subsection{Impact of a particle-binder torus against a rigid wall}

A particle-binder torus, that is made of a circular chain of $N_p$ spherical particles with diameter $D_p$ and that has an inner diameter of $D_t-D_p$ and an outer diameter of $D_t+D_p$, is initiated with an initial velocity field $v(t=0) = [-v_0, 0, 0]$ against a rigid wall with normal $[1,0,0]$ located at $x_w=0$. The diameter of the torus is $D_t = 3$, and it is comprised of eighty particles ($N_p=80$) of $m_{(i)}=1$, and $J_{(i)}=1$ which are in point contact in the undeformed configuration, that is  $D_p=D_{(i)}=D_t \sin(\pi/N_p)$. The stiffness coefficients of the binder are $K^\mathrm{m}_{(ij)}=10$, $K^\mathrm{a}_{(ij)}= K^\mathrm{s}_{(ij)} = 200$, and the particle stiffness coefficient is $K^\mathrm{pp}_{(ij)}=2,100$. The initial velocity is $v_0=1$.  Figure~\ref{Fig-TorusImpact-XY} shows a sequence of deformed states of the solid torus experiencing large deformation and vibrations after impact, for $h=10^{-3}$ and $t_N=35$. 

\begin{figure}[H]
	\centering{
			\includegraphics[scale=0.66]{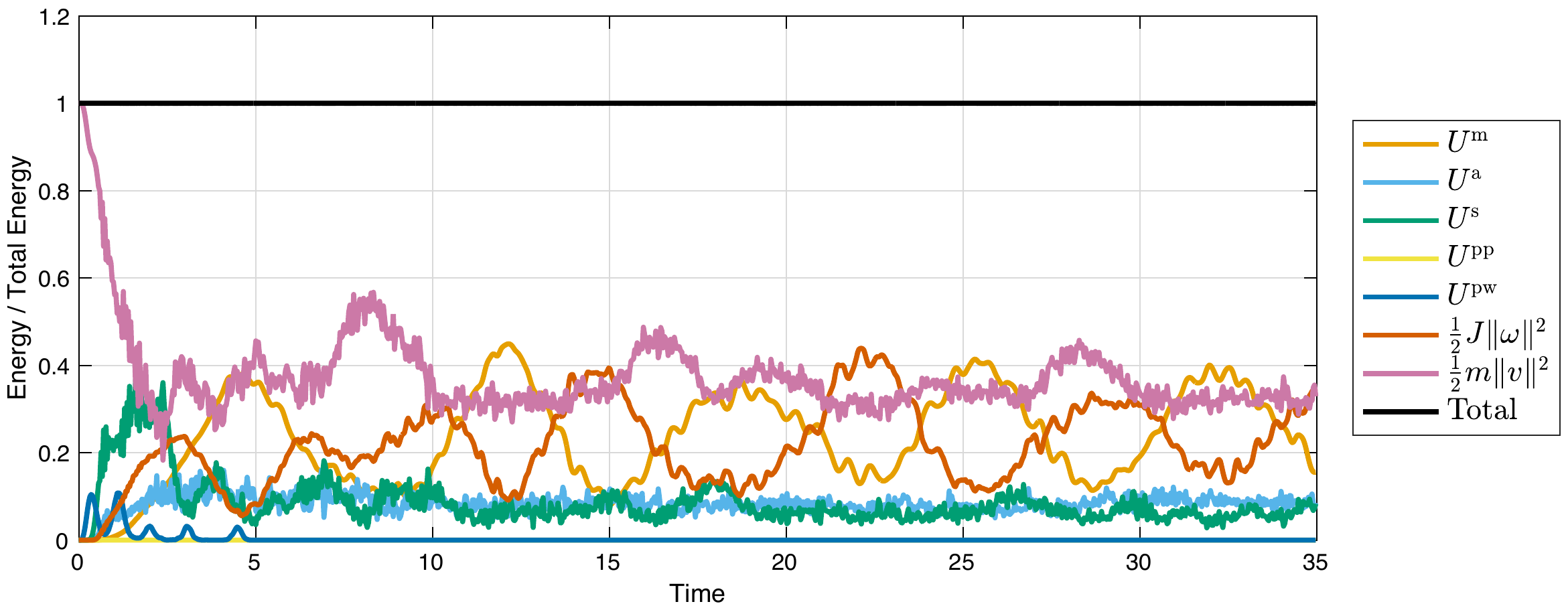}
	}
	\caption{For $h=10^{-3}$ and $t_N=25$,  the second-order, explicit time integrator given by discrete map \eqref{Eqn-VTI-2} exhibits excellent near-conservation of total energy over exponentially long times.}
	\label{Fig-TorusImpact-Energy-1}
\end{figure}

\begin{figure}[t]
	\centering{
		\includegraphics[scale=0.66]{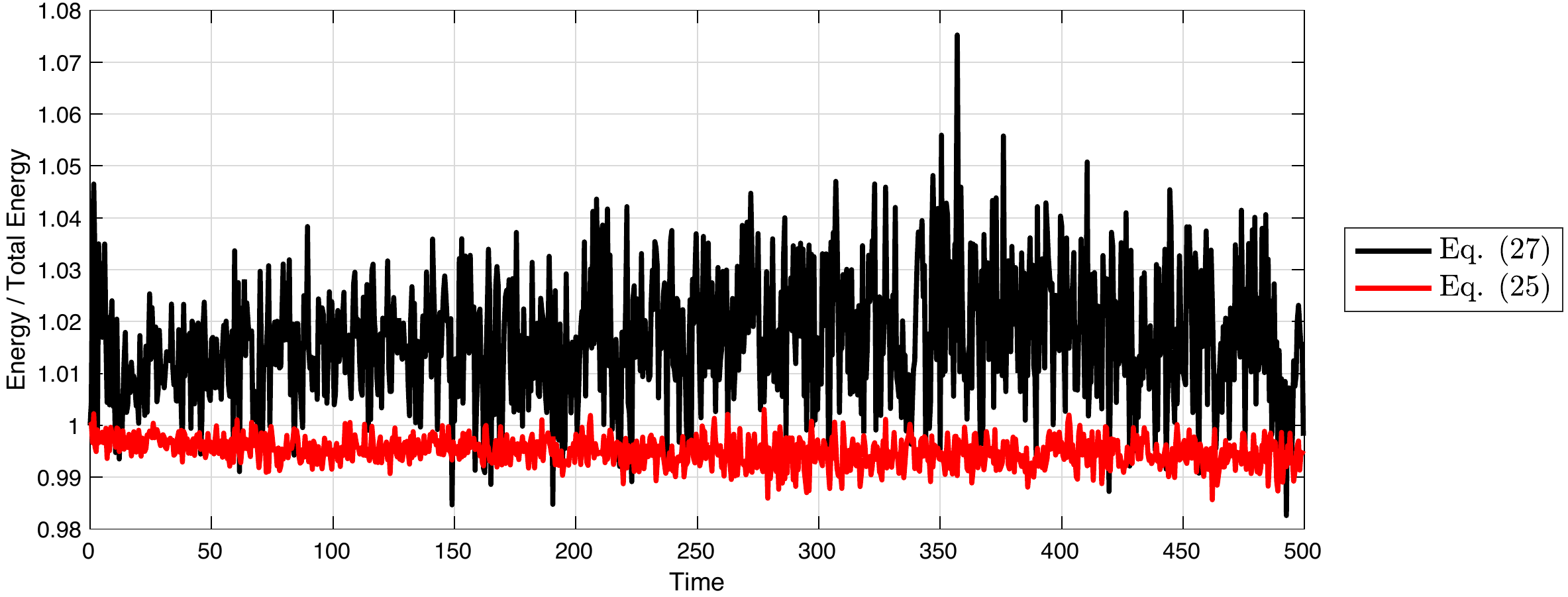}
	}
	\caption{Excellent near-conservation of total energy over exponentially long times exhibited by first-order \eqref{Eqn-VTI-3} and second-order \eqref{Eqn-VTI-2} variational time integrators, for $h=10^{-2}$ and $t_N=500$.}
	\label{Fig-TorusImpact-Energy-2}
\end{figure}

Figure \ref{Fig-TorusImpact-Energy-1} shows the excellent near-conservation properties of the explicit time integrator given discrete map \eqref{Eqn-VTI-2}. It is evident that, upon impact of the deformable torus against the wall, the initial translational kinetic energy ($m \|v\|^2/2$) is transfer to rotational kinetic energy ($J\|w\|^2/2$) and vibrational energy in the form of bending ($U^\mathrm{m}$), shear ($U^\mathrm{s}$), and axial ($U^\mathrm{a}$) interactions, as well as strain energies between particles ($U^\mathrm{pp}$) and between particles and the rigid wall ($U^\mathrm{pw}$). Despite the delicate balance between different energy forms during the dynamic response, the time integrator exhibits excellent stability over exponentially long times (see Figure \ref{Fig-TorusImpact-Energy-2}). This behavior is shared by all other variational time integrators derived from a discrete Lagrangian function.

\begin{figure}[H]
	\centering{
	\begin{tabular}{cc}
		\includegraphics[clip, trim=0.9cm 6.6cm 1.7cm 6.5cm, scale=0.41]{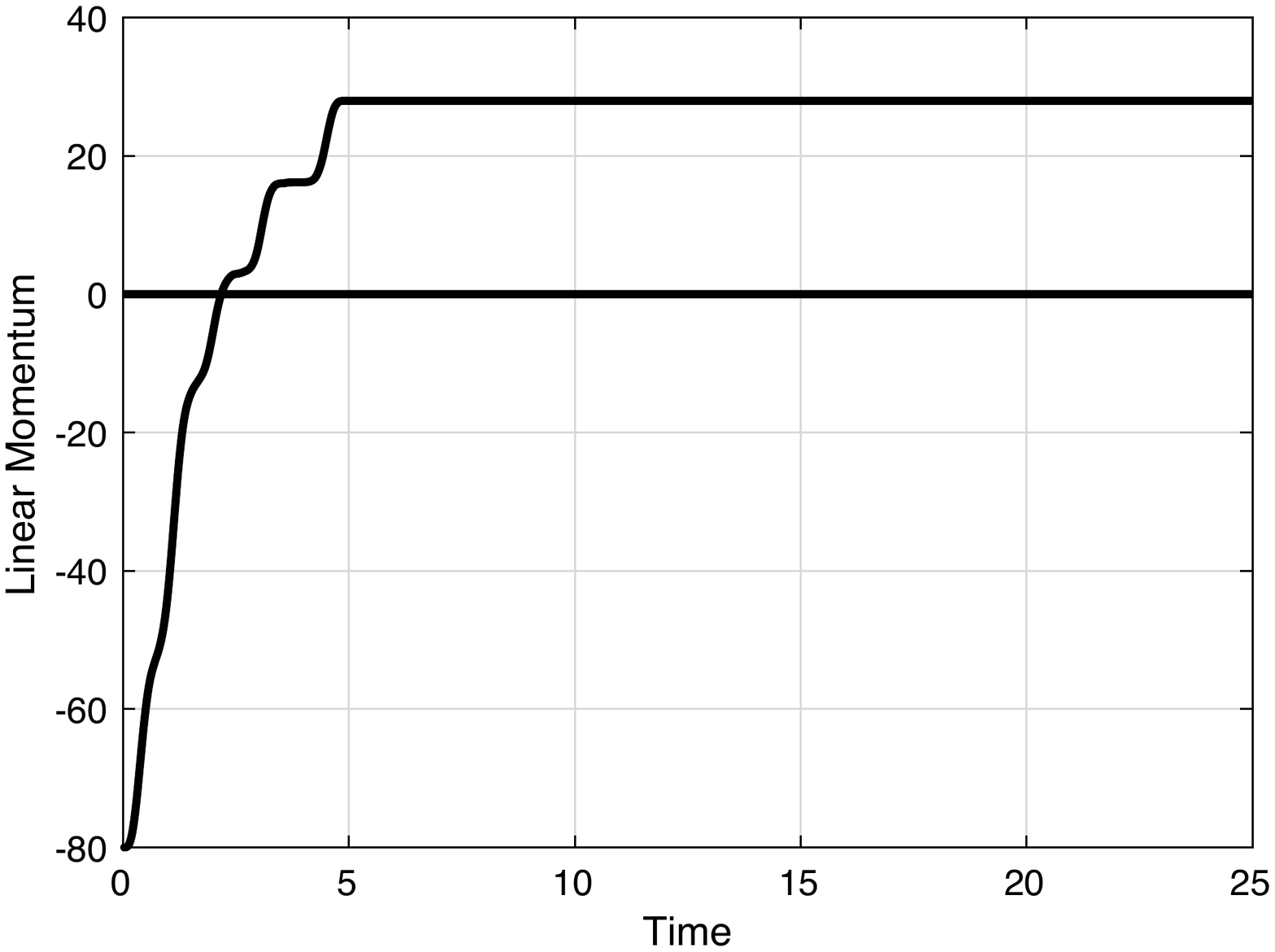}
		&
		\includegraphics[clip, trim=0.9cm 6.6cm 1.7cm 6.5cm, scale=0.41]{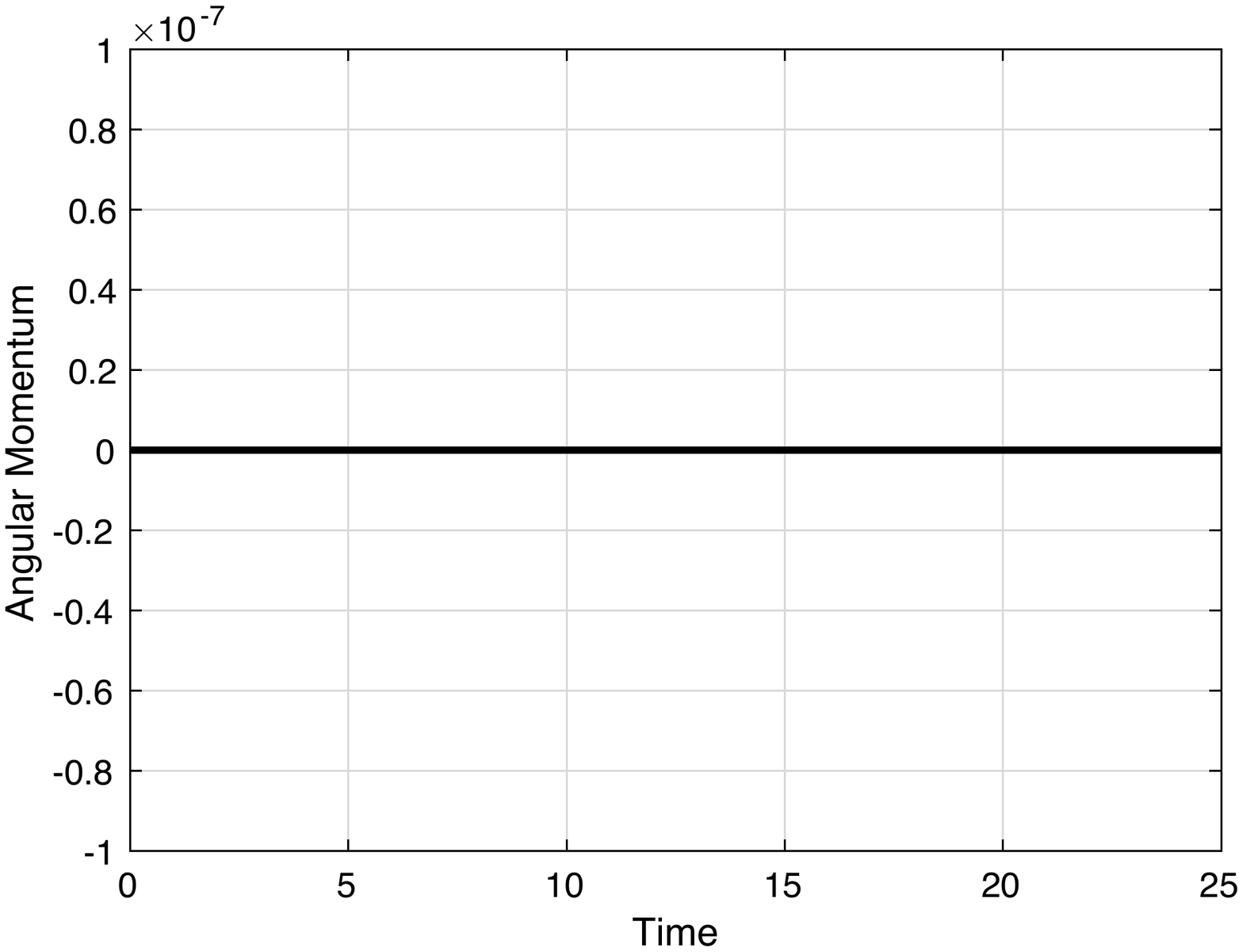}
	\end{tabular}
	}
	\caption{For $h=10^{-3}$ and $t_N=25$, both explicit time integrators exhibit exact conservation, to machine precision, of linear and angular momenta.}
	\label{Fig-TorusImpact-Momenta}
\end{figure}

Figure \ref{Fig-TorusImpact-Momenta} shows that the second-order, explicit integrators given by discrete maps \eqref{Eqn-VTI-1} and \eqref{Eqn-VTI-2} exhibit exact conservation, to machine precision, of linear and angular momenta. This behavior is shared by all other discrete momentum maps derived using a discrete Legendre transformation on the discrete Lagrangian that generates the variational time integrator.

\begin{figure}[H]
	\centering{
		\begin{tabular}{cc}
			\includegraphics[clip, scale=0.40]{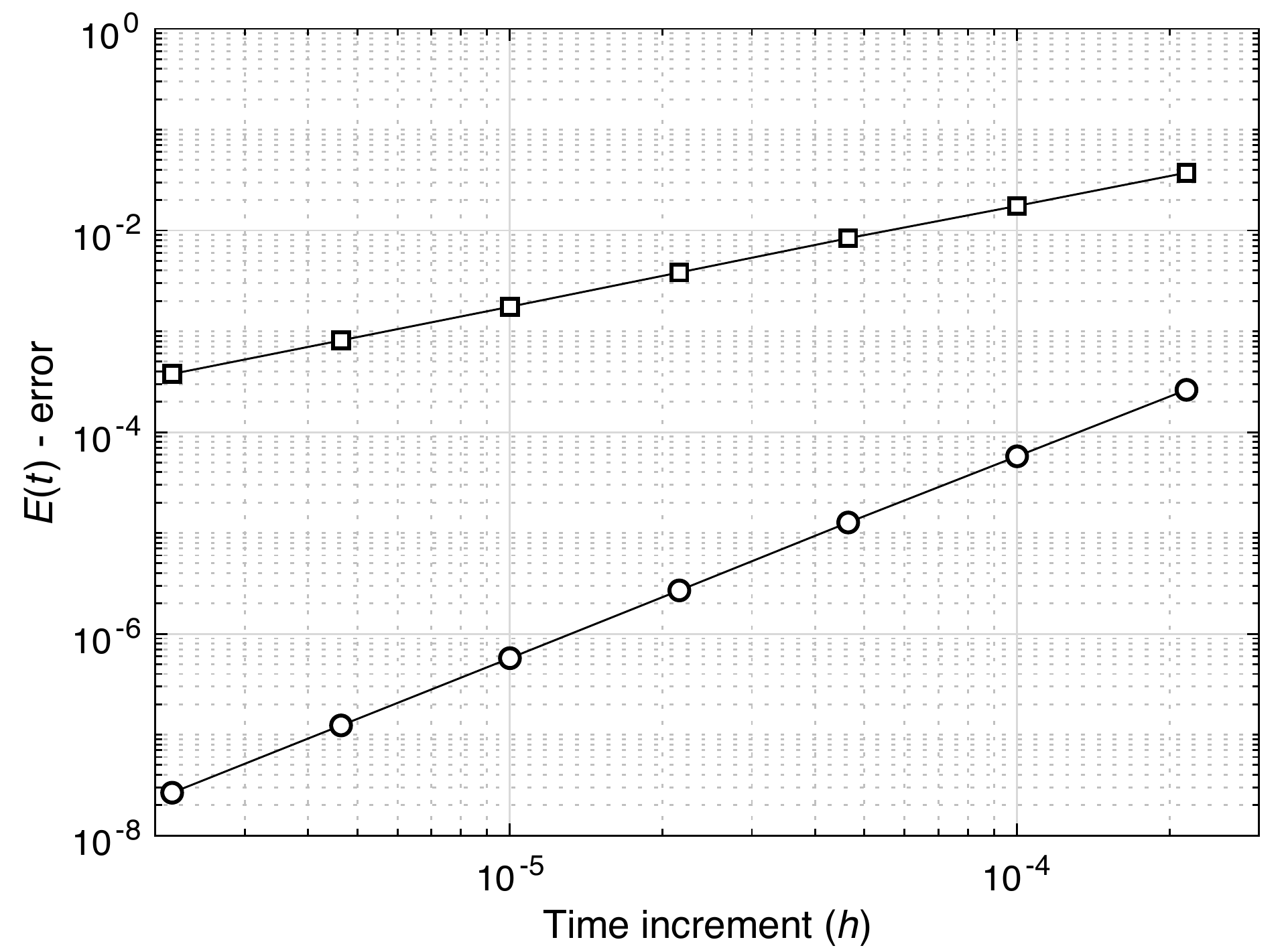}
			&
			\includegraphics[clip, scale=0.40]{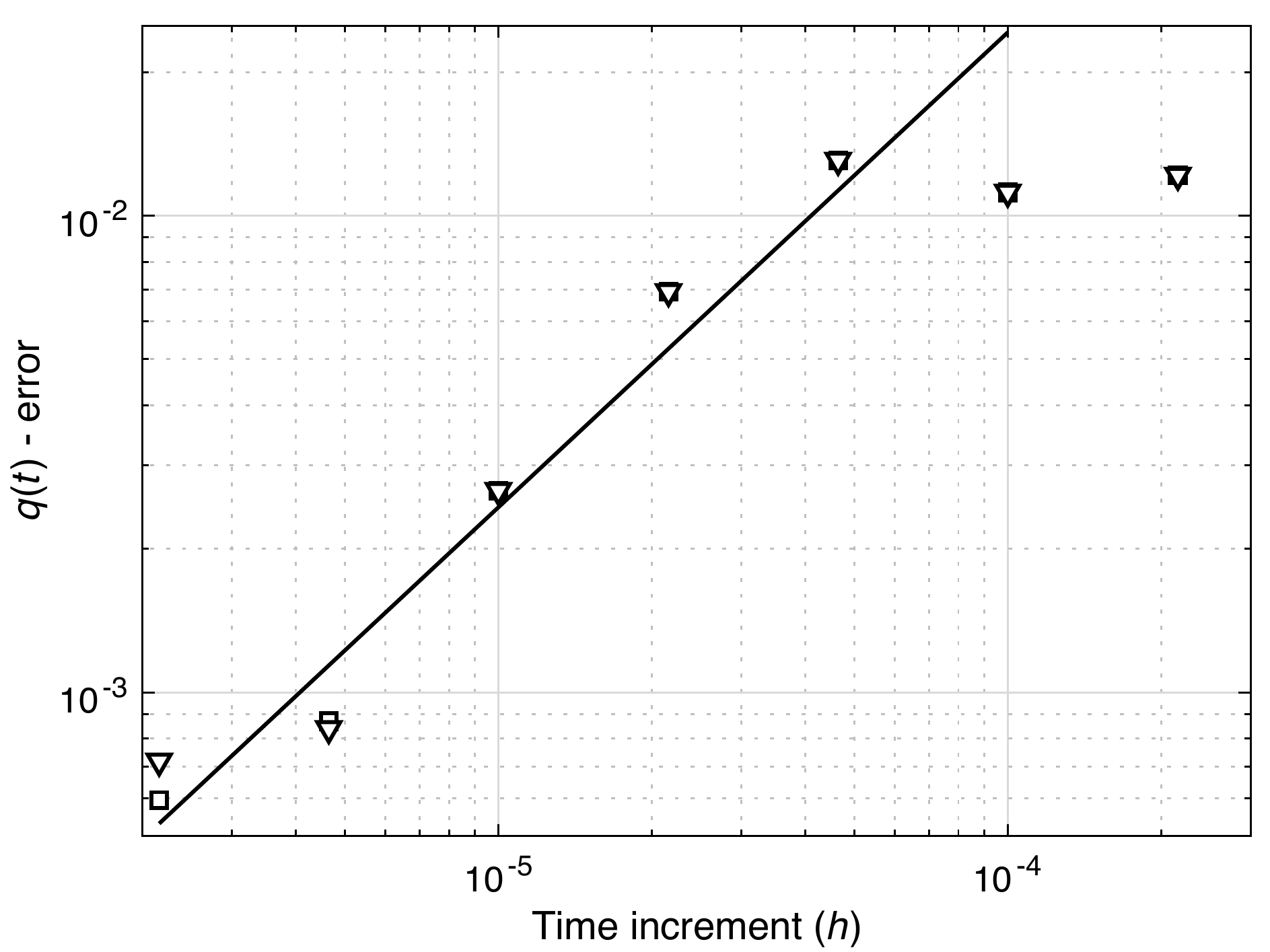}
		\end{tabular}
	}
	\caption{For $t_N=10$, both second explicit time integrators---discrete map \eqref{Eqn-VTI-1} in $\Circle$, \eqref{Eqn-VTI-2} in $\triangledown$---exhibit (i) quadratic convergence of energy in the $H^0(0,10)$-norm using the  $E(t)\text{-error}$, and (ii) linear convergence of trajectories in the $H^1(0,10)$-norm using the $q(t)\text{-error}$. The first-order explicit time integrator---discrete map \eqref{Eqn-VTI-3} in $\square$---exhibits (i) linear convergence of energy in the $H^0(0,10)$-norm using the  $E(t)\text{-error}$, and (ii) linear convergence of trajectories in the $H^1(0,10)$-norm using the $q(t)\text{-error}$.}
	\label{Fig-TorusImpact-Convergence}
\end{figure}

Finally, Figure \ref{Fig-TorusImpact-Convergence} shows that the second-order, explicit integrators given by discrete maps \eqref{Eqn-VTI-1} and \eqref{Eqn-VTI-2} exhibit (i) quadratic convergence of energy in the $H^0$-norm using the $E(t)\text{-error}$, and (ii) linear convergence of trajectories in the $H^1$-norm using the $q(t)\text{-error}$. The figure also shows results for the first-order variational time integrator given by the discrete map \eqref{Eqn-VTI-3}.

\section{Summary and discussion}

By parametrizing the space of rotations using exponential local coordinates represented by a rescaled form of the Rodrigues rotation vector, we have proposed a discrete Lagrangian function to derived a variational time integrator in $\mathrm{SE}(3)$.  We restricted attention to (i) Lagrangian functions with a quadratic kinetic energy and potential energies that solely depend on positions and attitudes, and (ii) spherical bodies, with the origin of the body-fixed frame located at the center of mass. The resulting discrete Euler-Lagrange equations can be solved analytically and, thus, we have systematically derived an explicit, closed-form expression of the corresponding discrete Lagrangian map, i.e., of the variational time integrator. 

We have shown that the proposed second-order, explicit time integrator based on rescaled Rodrigues parameters preserves the momenta of the continuous dynamics, such as linear and angular momenta, and that it exhibits near-conservation of total energy over exponentially long times. By adopting a inertially fixed frame and retaining terms in discrete Lagrangian map up to second order in time, we have obtained a discrete map that exhibits the same mathematical structure of the explicit Newmark time integrator ($\beta=0$ and $\gamma=1/2$) or the velocity Verlet algorithm, both known to be variational time integrators. Furthermore, we have demonstrated, by numerical experimentation, that the preserving properties of the time integrator are not diminished by this approximation, but they are rather indistinguishable from each other. Specifically, these properties are born out by two examples, namely the dynamics on $\mathrm{SO}(3)$ of a three-dimensional pendulum, and the nonlinear dynamics on $\mathrm{SE}(3)^n$ that results from the impact of a particle-binder torus against a rigid wall. These examples showcase that the proposed integrator exhibits linear convergence of trajectories in the $H^1$-norm, i.e., in the Sobolev space $H^1:=W^{1,2}$, and quadratic convergence of energy in the $H^0$-norm.

We close by pointing out some limitations of our analysis and possible avenues for extensions of the approach.

First, we have restricted attention to spherical bodies. The  extension of the variational time integrator based on rescaled Rodrigues parameters to non-spherical particles and, therefore, the systematic investigation of closed-form analytical solutions of the relative angular position $\Delta \alpha_{k}$ for a general moment of inertia are worthwhile directions of future research.

Second, we have derived a time-stepping scheme that exhibits exact-conservation of the momenta but only near-conservation of total energy. It is suggestive to note that the close-form expressions derived using rescaled Rodrigues parameters may also enable the formulation of energy-stepping \cite{EnergySteppingIntegrators}, or force-stepping \cite{ForceSteppingIntegrators}, integrators in $\mathrm{SE}(3)^n$.  Discrete Lagrangian maps of an energy-stepping scheme are the exact trajectories of an approximate Lagrangians and, therefore, these trajectories preserve total energy and all momentum maps whose associated symmetries are preserved by the approximate Lagrangians.

\section*{Acknowledgments}

The authors gratefully acknowledge the support received from the U.S. Air Force Office of Scientific Research through Award No. FA9550-15-1-0102. Caroline Baker also acknowledges the Ross Fellowship and the Mechanical Engineering Ward A. Lambert Graduate Teaching Fellowship from Purdue University.

\bibliographystyle{plainnat}
\bibliography{all-biblatex}

\begin{thebibliography}{23}
\providecommand{\natexlab}[1]{#1}
\providecommand{\url}[1]{\texttt{#1}}
\expandafter\ifx\csname urlstyle\endcsname\relax
  \providecommand{\doi}[1]{doi: #1}\else
  \providecommand{\doi}{doi: \begingroup \urlstyle{rm}\Url}\fi

\bibitem[Agarwal and Gonzalez(2020)]{agarwal2020effects}
Ankit Agarwal and Marcial Gonzalez.
\newblock Effects of cyclic loading and time-recovery on microstructure and
  mechanical properties of particle-binder composites.
\newblock \emph{Journal of Applied Mechanics}, 87\penalty0 (10):\penalty0
  101008, 2020.

\bibitem[Agarwal and Gonzalez(2021)]{agarwal2021constitutive}
Ankit Agarwal and Marcial Gonzalez.
\newblock A finite-deformation constitutive model of particle-binder composites
  incorporating yield-surface-free plasticity.
\newblock \emph{Journal of Applied Mechanics}, in press, 2021.

\bibitem[Bake and Gonzalez(2021)]{Baker2021}
Caroline Bake and Marcial Gonzalez.
\newblock Particle mechanics approach to microstructural modeling of impact
  behavior in bonded granular solids.
\newblock \emph{Computational Particle Mechanics}, under review, 2021.

\bibitem[Campello et~al.(2011)Campello, Pimenta, and Wriggers]{Campello2011}
Eduardo de Morais~Barreto Campello, Paulo de~Mattos Pimenta, and P~Wriggers.
\newblock {An exact conserving algorithm for nonlinear dynamics with rotational
  DOFs and general hyperelasticity. Part 2: Shells}.
\newblock \emph{Computational Mechanics}, 48\penalty0 (2):\penalty0 195--211,
  2011.

\bibitem[Campello(2015)]{Campello2015}
Eduardo M.~B. Campello.
\newblock {A description of rotations for DEM models of particle systems}.
\newblock \emph{Computational Particle Mechanics}, 2\penalty0 (2):\penalty0
  109--125, 2015.

\bibitem[Gonzalez et~al.(2010{\natexlab{a}})Gonzalez, Schmidt, and
  Ortiz]{EnergySteppingIntegrators}
Marcial Gonzalez, Bernd Schmidt, and Michael Ortiz.
\newblock {Energy-stepping integrators in Lagrangian mechanics}.
\newblock \emph{International Journal for Numerical Methods in Engineering},
  82\penalty0 (2):\penalty0 205--241, 2010{\natexlab{a}}.

\bibitem[Gonzalez et~al.(2010{\natexlab{b}})Gonzalez, Schmidt, and
  Ortiz]{ForceSteppingIntegrators}
Marcial Gonzalez, Bernd Schmidt, and Michael Ortiz.
\newblock {Force-stepping integrators in Lagrangian mechanics}.
\newblock \emph{International journal for numerical methods in engineering},
  84\penalty0 (12):\penalty0 1407--1450, 2010{\natexlab{b}}.

\bibitem[Hairer et~al.(2006)Hairer, Lubich, and Wanner]{Hairer2006Geometric}
Ernst Hairer, Christian Lubich, and Gerhard Wanner.
\newblock \emph{{Geometric numerical integration, volume 31 of Springer Series
  in Computational Mathematics}}.
\newblock Springer-Verlag, Berlin, 2006.

\bibitem[Huynh(2009)]{Huynh2009Metrics}
Du~Q Huynh.
\newblock {Metrics for 3D rotations: Comparison and analysis}.
\newblock \emph{Journal of Mathematical Imaging and Vision}, 35\penalty0
  (2):\penalty0 155--164, 2009.

\bibitem[Kane et~al.(2000)Kane, Marsden, Ortiz, and
  West]{Kane2000NewmarkAlgorithm}
Couro Kane, Jerrold~E. Marsden, Michael Ortiz, and Matthew West.
\newblock {Variational integrators and the Newmark algorithm for conservative
  and dissipative mechanical systems}.
\newblock \emph{International Journal for Numerical Methods in Engineering},
  49\penalty0 (10):\penalty0 1295--1325, 2000.

\bibitem[Lee et~al.(2007)Lee, Leok, and McClamroch]{Lee2007}
Taeyoung Lee, Melvin Leok, and N~Harris McClamroch.
\newblock Lie group variational integrators for the full body problem in
  orbital mechanics.
\newblock \emph{Celestial Mechanics and Dynamical Astronomy}, 98\penalty0
  (2):\penalty0 121--144, 2007.

\bibitem[Leok and Shingel(2012)]{Leok2012}
Melvin Leok and Tatiana Shingel.
\newblock {General techniques for constructing variational integrators}.
\newblock \emph{Frontiers of Mathematics in China}, 2012.

\bibitem[Leok and Zhang(2011)]{Leok2011}
Melvin Leok and Jingjing Zhang.
\newblock {Discrete Hamiltonian variational integrators}.
\newblock \emph{IMA Journal of Numerical Analysis}, 2011.

\bibitem[Lew et~al.(2004)Lew, Marsden, Ortiz, and West]{Lew2004}
Adrian Lew, Jerrold~E. Marsden, Michael Ortiz, and Matthew West.
\newblock {An overview of variational integrators}.
\newblock \emph{Finite Element Methods: 1970's and Beyond}, pages 98--115,
  2004.

\bibitem[Marsden and West(2001)]{Marsden2001}
Jerrold~E. Marsden and Matthew West.
\newblock {Discrete mechanics and variational integrators}.
\newblock \emph{Acta Numerica 2001}, 10\penalty0 (2001):\penalty0 357--514,
  2001.

\bibitem[Pimenta et~al.(2008)Pimenta, Campello, and Wriggers]{Pimenta2008}
Paulo de~Mattos Pimenta, Eduardo de Morais~Barreto Campello, and P~Wriggers.
\newblock {An exact conserving algorithm for nonlinear dynamics with rotational
  DOFs and general hyperelasticity. Part 1: Rods}.
\newblock \emph{Computational Mechanics}, 42\penalty0 (5):\penalty0 715--732,
  2008.

\bibitem[Poorsolhjouy and Gonzalez(2021)]{Poorsolhjouy2021}
Payam Poorsolhjouy and Marcial Gonzalez.
\newblock Particle mechanics approach to microstructural modeling of impact
  behavior in bonded granular solids.
\newblock \emph{Composite Structures}, in press, 2021.

\bibitem[Schaub and Junkins(2003)]{Schaub2003}
Hanspeter Schaub and John~L Junkins.
\newblock \emph{Analytical mechanics of space systems}.
\newblock AIAA, 2003.

\bibitem[Si et~al.(2002)Si, Little, and Lytton]{Si2002}
Zhiming Si, DN~Little, and RL~Lytton.
\newblock Characterization of microdamage and healing of asphalt concrete
  mixtures.
\newblock \emph{Journal of materials in civil engineering}, 14\penalty0
  (6):\penalty0 461--470, 2002.

\bibitem[Wang(2009)]{Wang2009}
Yucang Wang.
\newblock {A new algorithm to model the dynamics of 3-D bonded rigid bodies
  with rotations}.
\newblock \emph{Acta Geotechnica}, 4\penalty0 (2):\penalty0 117--127, 2009.

\bibitem[Wendlandt and Marsden(1997)]{Wendlandt1997}
Jeffrey~M. Wendlandt and Jerrold~E. Marsden.
\newblock {Mechanical integrators derived from a discrete variational
  principle}.
\newblock \emph{Physica D: Nonlinear Phenomena}, 106\penalty0 (3-4):\penalty0
  223--246, 1997.

\bibitem[Zhong and Marsden(1988)]{GeMarsden1988LiePoissonIntegrators}
Ge~Zhong and Jerrold~E. Marsden.
\newblock {Lie-Poisson Hamilton-Jacobi theory and Lie-Poisson integrators}.
\newblock \emph{Physics Letters A}, 133\penalty0 (3):\penalty0 134--139, 1988.

\bibitem[Zhu et~al.(1996)Zhu, Chang, and Rish]{Zhu1996}
Han Zhu, Ching~S. Chang, and Jeff~W. Rish.
\newblock {Normal and tangential compliance for conforming binder contact I:
  Elastic binder}.
\newblock \emph{International Journal of Solids and Structures}, 33\penalty0
  (29):\penalty0 4337--4349, 1996.

\end{thebibliography}

\end{document}